\documentclass[twocolumn,english,journal]{IEEEtran}
\usepackage[T1]{fontenc}
\usepackage{babel}
\usepackage{amsmath}
\usepackage{amssymb}
\usepackage{graphicx}
\usepackage{esint}
\usepackage[unicode=true,
 bookmarks=true,bookmarksnumbered=true,bookmarksopen=true,bookmarksopenlevel=1,
 breaklinks=false,pdfborder={0 0 0},pdfborderstyle={},backref=false,colorlinks=false]
 {hyperref}
\hypersetup{pdftitle={Your Title},
 pdfauthor={Your Name},
 pdfpagelayout=OneColumn, pdfnewwindow=true, pdfstartview=XYZ, plainpages=false}

\makeatletter
\newtheorem{assumption}{Assumption}
\newtheorem{problem}{Problem}
\newtheorem{thm}{Theorem}
\newtheorem{lemma}{Lemma}
\newtheorem{definitn}{Definition}
\newtheorem{prop}{Proposition}

\usepackage[caption=false,font=footnotesize]{subfig}

\makeatother

\begin{document}
\title{Online Optimal State Feedback Control of Linear Systems over Wireless
MIMO Fading Channels}
\author{Songfu Cai, \textsl{Member, IEEE}, Vincent K. N. Lau, \textit{Fellow,
IEEE}\\
Department of Electronic and Computer Engineering \\
The Hong Kong University of Science and Technology\\
Clear Water Bay, Kowloon, Hong Kong \\
Email: \{eesfcai, eeknlau\}@ust.hk}
\maketitle
\begin{abstract}
We consider the optimal control of linear systems over wireless MIMO
fading channels, where the MIMO wireless fading and random access
of the remote controller may cause intermittent controllability or
uncontrollability of the closed-loop control system. We formulate
the optimal control design over random access MIMO fading channels
as an infinite horizon average cost Markov decision process (MDP),
and we propose a novel state reduction technique such that the optimality
condition is transformed into a time-invariant reduced-state Bellman
optimality equation. We provide the closed-form characterizations
on the existence and uniqueness of the optimal control solution via
analyzing the reduced-state Bellman optimality equation. Specifically,
in the case that the closed-loop system is almost surely controllable,
we show that the optimal control solution always exists and is unique.
In the case that MIMO fading channels and the random access of the
remote controller destroy the closed-loop controllability, we propose
a novel controllable and uncontrollable positive semidefinite (PSD)
cone decomposition induced by the singular value decomposition (SVD)
of the MIMO fading channel contaminated control input matrix. Based
on the decomposed fine-grained reduced-state Bellman optimality equation,
we further propose a closed-form sufficient condition for both the
existence and the uniqueness of the optimal control solution. The
closed-form sufficient condition reveals the fact that the optimal
control action may still exist even if the closed-loop system suffers
from intermittent controllability or almost sure uncontrollability.
We further propose a novel stochastic approximation (SA)-based online
learning algorithm that can learn the optimal control action on the
fly based on the plant state observations. We derived a closed-form
sufficient condition that guarantees the almost sure convergence of
the online learning algorithm to the optimal control solution. The
proposed scheme is also compared with various baselines, and we show
that significant performance gains can be achieved.
\end{abstract}
\begin{IEEEkeywords}
Online learning, optimal control, wireless MIMO fading channels, Markov
decision process, uncontrollable linear systems, Lyapunov stability
analysis.
\end{IEEEkeywords}

\section{Introduction}

\begin{figure}[t!]
\centering{}\includegraphics[width=1\columnwidth]{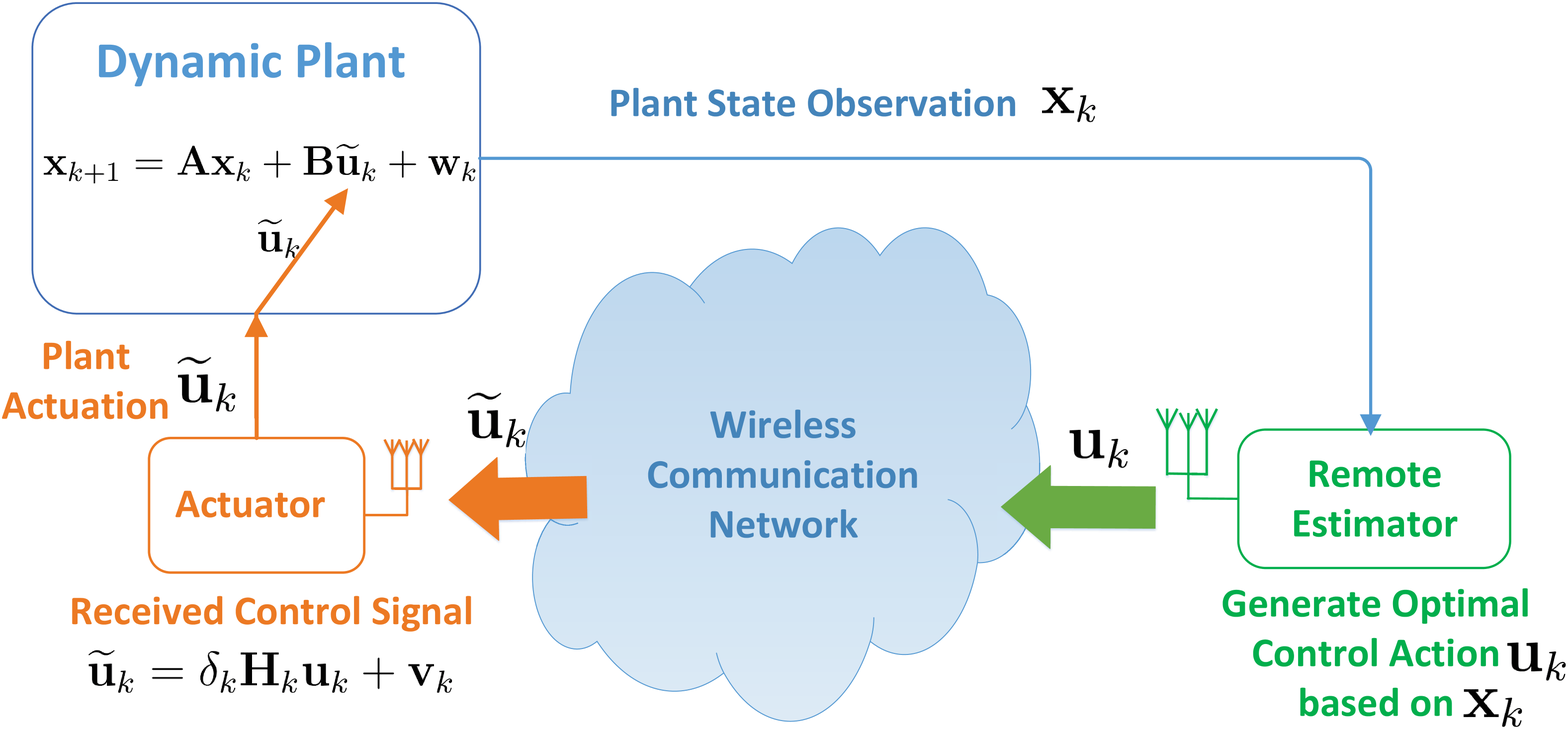} \caption{\label{fig-system-topo} Illustration of a closed-loop state feedback
control system, where the remote controller transmits its control
action to the actuator over a MIMO wireless communication network.}
\end{figure}

Optimal control has received considerable attention in both academia
and industry in recent years. A wide spectrum of applications of optimal
control can be found in areas such as aerospace control, industrial
flotation process control, automated vehicle systems, and robotics
and manufacturing systems {[}1{]}-{[}3{]}. A typical closed-loop feedback
control system consists of a dynamic plant (with potentially unstable
dynamics), a remote controller, and an actuator, as illustrated in
Fig. \ref{fig-system-topo}. Specifically, the remote controller generates
the real-time plant control action based on the instantaneous plant
state observation. The remote controller then transmits the plant
control action to the actuator over a wireless communication network.
The actuator, which is collocated with the plant, applies its received
control signals for plant actuation. The wireless network in-between
the remote controller and the actuator will have significant impacts
on the closed-loop control performance because it introduces various
degradations, such as wireless fading, packet errors and latency {[}4{]}-{[}6{]}.
As a result, it is important to incorporate the impairments in the
wireless networks into the optimal control design at the remote controller.

Recently, there have been some works on optimal control over static
channels. Specifically, in {[}7{]}, the authors consider the optimal
real-time control of wind turbines. The optimal control action is
obtained numerically using dynamic programming by maximizing the wind
energy capture. In {[}8{]}, the authors consider the problem of optimal
control with transmission power management. Exploiting the information
structure at the controller, the authors show that the linear quadratic
regulator (LQR) control law is optimal. In {[}9{]}, the authors consider
optimal control design by iteratively solving the associated Bellman
optimality equation. Specifically, under an initial stabilizing control
policy, the optimal control action is obtained via policy iteration,
whereas the initial stabilizing control policy is obtained by solving
a linear matrix inequality (LMI). In {[}10{]}, the authors propose
an adaptive optimal control design using adaptive dynamic programming
such that the a priori knowledge of an initial stabilizing control
policy is no longer required. In {[}11{]}, the authors propose a mixed
mode value and policy iteration algorithm to obtain the optimal control
action, which also avoids the need for an initial stabilizing control.
However, in all the aforementioned works {[}7{]}-{[}11{]}, the associated
optimal control gain is static due to the consideration of the static
communication channel between the remote controller and the actuator.
Therefore, these existing approaches cannot achieve optimal control
over dynamic wireless fading channels and brute-force applications
of the existing methods {[}7{]}--{[}11{]} may even cause plant instability.

There are also several works considering the optimal control over
wireless communication channels. In {[}12{]}, the authors model the
wireless communication channel between the remote controller and the
actuator as a point-to-point packet dropping channel with a fixed
number of packet drops in a certain time interval. The authors provide
closed-form characterizations on the optimal control action that minimizes
a finite-horizon LQR cost. In {[}13{]} and {[}14{]}, the authors consider
the finite burst of consecutive packet dropouts, and the optimal control
action is obtained using the potential learning approach. In {[}15{]},
the authors considers the i.i.d. random Bernoulli packet loss channels.
By using a separation principle, the authors show that the optimal
controller is a standard LQR controller. However, such a packet-dropping
channel model in {[}12{]}-{[}15{]} is an oversimplification of the
impairments introduced in practical wireless MIMO fading channels.
In {[}16{]}, the authors model the time varying wireless communication
channel between the remote controller and the actuator as a switched
linear system (SLS) with a finite number of switching states. However,
in practice, the wireless MIMO fading channel is switching in a continuous
state space with uncountably many realizations, where the analysis
approaches of SLS are thus not applicable. In {[}17{]} and {[}18{]},
the authors consider the optimal control over continuous fading channels
but only diagonal fading channels are considered. Moreover, the system
is required to be controllable at every timeslot, which is a very
restrictive requirement. Unfortunately, the MIMO fading channels are
far more complicated than the diagonal fading channels in the sense
that MIMO fading channels can destroy the system controllability and
the resultant system may not be controllable at every timeslot. The
impacts of intermittent uncontrollability on the closed-loop stability
and optimal control design have not been considered.

In this paper, we consider the optimal control over MIMO wireless
fading channels, where the wireless fading may cause intermittent
uncontrollability of the closed-loop control system. We propose a
novel online learning algorithm that converges almost surely to the
optimal control action. The following summarizes the key contributions
of the work.
\begin{itemize}
\item \textbf{Reduced State Bellman Equation: }We formulate the optimal
control design over MIMO fading channels as an infinite horizon average
cost Markov decision process (MDP). Due to the dynamic MIMO channel
fading, the state space of the stochastic optimal control problem
has to be extended to include both the channel state information (CSI)
and the plant state information (PSI). As a result, learning the value
function of the associated Bellman optimality equation is more challenging
due to the expanded state space. To overcome this obstacle, we propose
a novel\textsl{ $\emph{state reduction technique}$} so that the optimality
condition is transformed into a time-invariant\textsl{ }$\emph{reduced-state Bellman optimality equation}$.
\item \textbf{Existence and Uniqueness of Optimal Control with Intermittent
Controllability or Almost Sure Uncontrollability:} We provide closed
form characterization of the sufficient condition for the existence
and uniqueness of the optimal control solution over wireless MIMO
fading channels. We propose a novel controllable and uncontrollable
PSD cone decomposition technique. We show that the reduced state Bellman
equation has an unique solution if we have $\emph{almost sure controllability}$.
For the cases of $\emph{intermittent controllability}$ and $\emph{almost sure uncontrollability}$,
we provide closed-form sufficient condition for existence and uniqueness
of the optimal control solution.
\item \textbf{Online Learning the Optimal Control Action and the Convergence
Analysis:} We propose a novel SA-based online learning algorithm that
can learn the optimal control action $\mathbf{u}_{k}$ on the fly
based on the plant state observations $\mathbf{x}_{k}$. The convergence
of the proposed SA-based online learning algorithm is characterized
via analyzing the associated limiting ordinary differential equation
(ODE), where the MIMO fading channel states and the kernel of the
value function are tightly coupled together in a highly nonlinear
manner. To address this challenge, we introduce a virtual fixed-point
process, for which the state trajectory is arbitrarily close to the
trajectory of the limiting ODE. By analyzing the fixed-point operator
associated with the virtual fixed-point process, we derive a closed-form
sufficient condition for the convergence of the limiting ODE, which
in turn renders the almost sure convergence of the online learning
algorithm to the optimal control solution.
\end{itemize}

\textsl{Notation}: Uppercase and lowercase boldface denote matrices
and vectors, respectively. The operator $(\cdot)^{T}$, $(\cdot)^{\dagger}$,
$(\cdot)^{H}$, $\left\lfloor \cdot\right\rfloor $, $\mathrm{Tr\left(\cdot\right)}$,
and $\mathrm{Re}\left\{ \cdot\right\} $ is the transpose, element-wise
conjugate, conjugate transpose, floor function, trace, and real part,
respectively. $\mathbf{0}_{m\times n}$ denotes $m\times n$ dimensional
matrices with all the elements being zero. $||\mathbf{A}||$ denotes
the spectrum norm of matrix $\mathbf{A}$. $\rho\left(\mathbf{A}\right)$
denote the spectral radius of matrix $\mathbf{A}$. $||\mathbf{a}||$
denotes the $l_{2}$ norm of vector $\mathbf{a}$. $||\mathbf{f}||$
denotes the operator norm of the operator $\mathbf{f}\left(\cdot\right)$.
$\left(\mathbf{a}\right)_{i}$ denotes the $i$-th entry of vector
$\mathbf{a}$. $\mathbf{A}_{ij}$ denotes the element in the $i$-th
row and $j$-th column of matrix $\mathbf{A}$. $\left(\mathbf{A}\right)_{i}$
denotes the $i$-th order leading principal submatrix of $\mathbf{A}$.
$\mathbf{\left(A\right)}_{i:j;l:m}$ denotes the $(j-i+1)\times(m-l+1)$
dimensional block submatrix of $\mathbf{A}$ with the first element
being $\mathbf{A}_{ij}$. $\mathbb{S}_{+}^{S}$ denotes the set of
$S\times S$ dimensional positive definite matrices. $\mathbb{S}^{S}$
denotes the set of $S\times S$ dimensional positive semidefinite
matrices. $\mathbb{R}^{m\times n}$ ($\mathbb{C}^{m\times n}$) represents
the set of $m\times n$ dimensional real (complex) matrices.

\section{System Model}

In this section, we introduce the architecture and key components
of the closed-loop feedback control system and formulate the optimal
control problem.

\subsection{Dynamic Plant Model}

A typical closed-loop feedback control system is a geographically
distributed system, wherein a potentially unstable dynamic plant,
an actuator, and a remote controller are connected through a wireless
communication network, as illustrated in Fig. \ref{fig-system-topo}.
The dynamic plant is modeled as a linear dynamic system, which is
described by a set of first order coupled linear difference equations
representing the evolution of the state variables. The dynamic evolution
of the plant state $\mathbf{x}_{k}$ is summarized below.

\begin{assumption}
\emph{(Dynamic Plant Model)} The plant state $\mathbf{x}_{k}$ follows
the dynamic evolution of $\mathbf{x}_{k+1}=\mathbf{A}\mathbf{x}_{k}+\mathbf{B}\widetilde{\mathbf{u}}_{k}+\mathbf{w}_{k}$,
$\ k\geq0$, where $\mathbf{x}_{k}\in\mathbb{R}^{S\times1}$ is the
plant state process, $S$ is the plant state dimension, $\mathbf{x}_{0}$
is the initial state vector, $\widetilde{\mathbf{u}}_{k}\in\mathbb{R}^{N_{r}\times1}$
is the actuation control input signal, $\mathbf{A}\in\mathbb{R}^{K\times K}$,
$\mathbf{B}\in\mathbb{R}^{S\times N_{r}}$, and $\mathbf{w}_{k}\in\mathbb{R}^{S\times1}$
is the plant noise with zero mean and finite covariance matrix $\mathbf{W}$\footnote{There exist a bounded constant $W$ such that $\mathbf{W}\leq W\mathbf{I}$.}.
The plant state transition matrix $\mathbf{A}$ contains possibly
unstable eigenvalues.
\end{assumption}

\subsection{Wireless Communication Model}

We model the wireless communication channel between the multi-antenna
remote controller and the actuator as a wireless MIMO fading channel.
Using multiple-antenna techniques, the $N_{t}$- antenna controller
transmits its control action $\mathbf{u}_{k}$ to the $N_{r}$-antenna
actuator through spatial multiplexing. At the $k$-th time slot, the
received control signal $\widetilde{\mathbf{u}}_{k}\in\mathbb{R}^{N_{r}\times1}$
at the actuator is given by
\begin{align}
 & \widetilde{\mathbf{u}}_{k}=\delta_{k}\mathbf{H}_{k}\mathbf{u}_{k}+\mathbf{v}_{k},
\end{align}
where $\mathbf{H}_{k}\in\mathbb{R}^{N_{r}\times N_{t}}$ is the MIMO
channel fading matrix, $\delta_{k}\in\left\{ 0,\ 1\right\} $ is the
$random\ access\ variable$ indicating whether the remote controller
is active to transmit its control action $\mathbf{u}_{k}\in\mathbb{R}^{N_{t}\times1}$
or not, and $\mathbf{v}_{k}\sim\mathcal{N}\left(0,\mathbf{I}_{N_{r}}\right)$
is the additive Gaussian channel noise. We have the following assumption
on $\mathbf{H}_{k}$.

\begin{assumption}
\emph{(MIMO Wireless Fading Channel Model)} The random MIMO channel
realization $\mathbf{H}_{k}$ remains constant within each time slot
and is i.i.d. over slots. Each element of $\mathbf{H}_{k}$ is i.i.d.
Gaussian distributed with zero mean and unit variance.
\end{assumption}

\subsection{Optimal Control Formulation}

The optimal state feedback control for a linear time-invariant (LTI)
system has been widely studied in existing literature {[}19{]}-{[}23{]}.
Specifically, a control policy $\pi$ consists of a sequence of mappings
$\pi=\left\{ \mathrm{\Omega}^{0},\mathrm{\Omega}^{1},\cdots\right\} $.
The mapping $\mathrm{\Omega}^{k}:\mathbb{R}^{S\times1}\rightarrow\mathbb{R}^{N_{t}\times1}$
at the $k$-th timeslot is a mapping from the plant state $\mathbf{x}_{k}$
to the control action $\mathbf{u}_{k}$, i.e., $\mathbf{u}_{k}=\mathrm{\Omega}^{k}\left(\mathbf{x}_{k}\right)$.
In a noiseless plant case, the optimal control is formulated as an
infinite horizon LQR total cost minimization problem {[}19{]} and
{[}20{]}:
\begin{align}
\min_{\pi} & \ \ \ \mathcal{J}^{\pi}=\sum_{k=0}^{\infty}r\left(\mathbf{x}_{k},\mathbf{u}_{k}\right),\label{eq:std-formulation-1}\\
s.t. & \ \ \ \mathbf{x}_{k+1}=\mathbf{A}\mathbf{x}_{k}+\mathbf{B}\mathbf{u}_{k},\nonumber 
\end{align}
where $r\left(\mathbf{x}_{k},\mathbf{u}_{k}\right)=\mathbf{x}_{k}^{T}\mathbf{Q}\mathbf{x}_{k}+\mathbf{u}_{k}^{T}\mathbf{R}\mathbf{u}_{k}$
is the per-stage cost reflecting the quadratic cost of state $\mathbf{x}_{k}^{T}\mathbf{Q}\mathbf{x}_{k}$
and the control cost $\mathbf{u}_{k}^{T}\mathbf{R}\mathbf{u}_{k}$,
and $\mathbf{Q}\in\mathbb{S}_{+}^{S}$ and $\mathbf{R}\in\mathbb{S}_{+}^{N_{r}}$
are the weighting matrices.

On the other hand, in a noisy plant system case, the infinite sum
$\sum_{k=0}^{\infty}r\left(\mathbf{x}_{k},\mathbf{u}_{k}\right)$
becomes unbounded and is not well-defined. As a result, the infinite
horizon ergodic control formulation {[}21{]}-{[}23{]} has to be adopted.
In this case, the optimal control problem is formulated as an ergodic
cost minimization problem {[}21{]}-{[}23{]}:

\emph{
\begin{align}
\min_{\pi} & \ \ \ \mathcal{J}^{\pi}=\limsup_{K\rightarrow\infty}\frac{1}{K}\sum_{k=0}^{K}\mathbb{E}\left[r\left(\mathbf{x}_{k},\mathbf{u}_{k}\right)\right]\label{eq:std-formulation}\\
s.t. & \ \ \ \mathbf{x}_{k+1}=\mathbf{A}\mathbf{x}_{k}+\mathbf{B}\mathbf{u}_{k}+\mathbf{w}_{k}.\nonumber 
\end{align}
}

Note that in the above standard formulations (\ref{eq:std-formulation-1})
and (\ref{eq:std-formulation}), the system dynamics are required
to be linear and time-invariant. However, when the random access of
the remote controller $\delta_{k}$ and the MIMO fading channel $\mathbf{H}_{k}$
are considered, the equivalent plant state dynamics are given by
\begin{align}
\mathbf{x}_{k+1} & =\mathbf{A}\mathbf{x}_{k}+\delta_{k}\mathbf{B}\mathbf{H}_{k}\mathbf{u}_{k}+\mathbf{B}\mathbf{v}_{k}+\mathbf{w}_{k},
\end{align}
which is linear but time varying due to the controller random access
process $\left\{ \delta_{k},k\geq0\right\} $ and the random MIMO
channel fading process $\left\{ \mathbf{H}_{k},k\geq0\right\} $.
As a result, the existing formulations cannot be directly applied
to our case.

In order to formulate stochastic optimal control for the linear and
time varying (LTV) system, we first extend the system state from $\mathbf{x}_{k}\in\mathbb{R}^{S\times1}$
to $\mathbf{\mathbf{S}}_{k}=\left(\mathbf{x}_{k},\mathbf{H}_{k},\delta_{k}\right)\in\mathbb{R}^{S\times1}\times\mathbb{R}^{N_{r}\times N_{t}}\times\left\{ 0,\ 1\right\} $,
where the extended state $\mathbf{\mathbf{S}}_{k}$ incorporates the
PSI $\mathbf{x}_{k}$, the CSI $\mathbf{H}_{k}$, and the controller
random access state $\delta_{k}$. In this case, the control policy
$\pi$ is a sequence of mappings $\pi=\left\{ \mathrm{\Omega}^{0},\Omega^{1},\cdots\right\} $,
where the mapping $\Omega^{k}:\mathbb{R}^{S\times1}\times\mathbb{R}^{N_{r}\times N_{t}}\times\left\{ 0,\ 1\right\} \rightarrow\mathbb{R}^{N_{t}\times1}$
at the $k$-th timeslot is a mapping from the extended state $\mathbf{\mathbf{S}}_{k}=\left(\mathbf{x}_{k},\mathbf{H}_{k},\delta_{k}\right)$
to the control action $\mathbf{u}_{k}$, i.e., $\mathbf{u}_{k}=\mathrm{\Omega}^{k}\left(\mathbf{\mathbf{S}}_{k}\right)$.
This physically means that the control action should be adaptive to
the realizations of the plant state $\mathbf{x}_{k}$ (reflecting
the urgency of the control), the channel state $\mathbf{H}_{k}$ (revealing
the transmission opportunities in the wireless MIMO channel), and
the random access of the controller $\delta_{k}$ (indicating transmission
urgency). Furthermore, the per-stage cost will need to include the
state cost $\mathbf{x}_{k}^{T}\mathbf{Q}\mathbf{x}_{k}$, the control
cost $\mathbb{E}\left[\left.\widetilde{\mathbf{u}}_{k}^{T}\mathbf{M}\widetilde{\mathbf{u}}_{k}\right|\mathbf{\mathbf{S}}_{k}\right],\mathbf{M}\in\mathbb{S}_{+}^{N_{r}},$
and the transmission cost $\mathbf{u}_{k}^{T}\mathbf{R}\mathbf{u}_{k}$,
and is given by
\begin{align}
r\left(\mathbf{\mathbf{S}}_{k},\mathbf{u}_{k}\right) & =\mathbf{x}_{k}^{T}\mathbf{Q}\mathbf{x}_{k}+\mathbf{u}_{k}^{T}\mathbf{R}\mathbf{u}_{k}+\mathbb{E}\left[\left.\widetilde{\mathbf{u}}_{k}^{T}\mathbf{M}\widetilde{\mathbf{u}}_{k}\right|\mathbf{\mathbf{S}}_{k}\right].
\end{align}

In addition, the extended state sequence $\left\{ \mathbf{\mathbf{S}}_{k}\right\} $
is a controlled Markov process with the transition kernel given by
\begin{align}
 & \mathrm{Pr}\left[\left.\mathbf{\mathbf{S}}_{k+1}\right|\mathbf{\mathbf{S}}_{k},\mathbf{u}_{k}\right]=\mathrm{Pr}\left[\left.\mathbf{\mathbf{H}}_{k+1}\right|\mathbf{\mathbf{S}}_{k},\mathbf{u}_{k}\right]\cdot\mathrm{Pr}\left[\left.\mathbf{\mathbf{x}}_{k+1}\right|\mathbf{\mathbf{S}}_{k},\mathbf{u}_{k}\right]\nonumber \\
 & =\mathrm{Pr}\left[\left.\mathbf{\mathbf{H}}_{k+1}\right|\mathbf{\mathbf{S}}_{k},\mathbf{u}_{k}\right]\mathrm{Pr}\left[\left.\mathbf{\mathbf{x}}_{k+1}\right|\mathbf{\mathbf{x}}_{k},\mathbf{u}_{k}\right].
\end{align}

Therefore, the optimal control over the wireless MIMO fading channels
can be formulated as an infinite horizon ergodic control problem w.r.t.
the extended state $\mathbf{\mathbf{S}}_{k}$, which is summarized
in the following Problem 1.

\begin{problem}
\textsl{(Optimal State Feedback Control Problem over Wireless MIMO
Fading Channels)}\label{problem-1}\emph{
\begin{align}
\min_{\pi} & \ \ \ \mathcal{J}^{\pi}=\limsup_{K\rightarrow\infty}\frac{1}{K}\sum_{k=0}^{K}\mathbb{E}\left[r\left(\mathbf{\mathbf{S}}_{k},\mathbf{u}_{k}\right)\right]\label{eq: prps-formulation}\\
s.t. & \ \ \ \mathbf{x}_{k+1}=\mathbf{A}\mathbf{x}_{k}+\delta_{k}\mathbf{B}\mathbf{H}_{k}\mathbf{u}_{k}+\mathbf{B}\mathbf{v}_{k}+\mathbf{w}_{k}.\nonumber 
\end{align}
}
\end{problem}

\section{Optimality Condition}

In this section, we first introduce the reduced state Bellman optimality
equation, which serves as an optimality condition for solving Problem
1. Based on the reduced state Bellman optimality equation, we further
provide a sufficient condition for the existence and uniqueness of
the optimal control action.

\subsection{Reduced State Bellman Optimality Equation}

Since the optimal state feedback control Problem \ref{problem-1}
is an infinite horizon ergodic control problem, the optimality condition
is given by the standard \textsl{Bellman equation} {[}24{]}, which
is summarized below.

\begin{thm}
(\textsl{Standard Optimality Condition for Problem 1)} \label{Thm: std-Bellman-eq}If
there exists a pair of $\left(\theta,V\left(\mathbf{S}_{k}\right)\right)$
that solves the following \textsl{Bellman optimality equation:}
\begin{align}
\theta+V\left(\mathbf{S}_{k}\right)= & \min_{\mathbf{u}_{k}}\bigg[r\left(\mathbf{\mathbf{S}}_{k},\mathbf{u}_{k}\right)\nonumber \\
 & +\sum_{\mathbf{S}_{k+1}}\mathrm{Pr}\left[\left.\mathbf{\mathbf{S}}_{k+1}\right|\mathbf{\mathbf{S}}_{k},\mathbf{u}_{k}\right]V\left(\mathbf{S}_{k+1}\right)\bigg],\forall\mathbf{S}_{k},\label{eq: standard Bellman Eq}
\end{align}
then:
\end{thm}

\begin{itemize}
\item For all initial extended states $\mathbf{S}_{0}$, $\theta=\min_{\mathbf{u}_{k}}\limsup_{K\rightarrow\infty}\frac{1}{K}\sum_{k=0}^{K}\mathbb{E}\left[r\left(\mathbf{\mathbf{S}}_{k},\mathbf{u}_{k}\right)\right]$
is the optimal average cost for Problem \ref{problem-1}, which is
independent of any extended states $\mathbf{S}_{0}$.
\item $V\left(\mathbf{S}_{k}\right)$ is the optimal value function for
the extended state $\mathbf{S}_{k}$.
\item The optimal control policy for Problem \ref{problem-1} is given by
$\mathbf{u}_{k}^{*}$, which attains the minimum of the R.H.S. of
(\ref{eq: standard Bellman Eq}).
\end{itemize}
\begin{IEEEproof}
Please see Appendix \ref{subsec:Proof-of-Thm1-Thm2}.
\end{IEEEproof}

There are various standard techniques such as value iteration {[}25{]}
and {[}26{]} or Q-learning {[}27{]}-{[}29{]} that can be used to solve
the \textsl{Bellman optimality equation} (\ref{eq: standard Bellman Eq}).
However, there are two challenges to solve. Firstly, there is a curse
of dimensionality in the extended state space $\left\{ \mathbf{S}_{k}\right\} $.
Specifically, the total dimension of the extended state $\mathbf{S}_{k}$
is $\left(1+S+N_{r}\cdot N_{t}\right)$, which can be huge when the
number of receive antennas at the actuator $N_{r}$ and the number
of transmit antennas at the remote controller $N_{t}$ are large.
As a result, if we brute-force applying the standard value iteration
to learn the value function $V\left(\mathbf{S}_{k}\right)$ {[}25{]}
and {[}26{]}, the domain of the value function contains $\left(1+S+N_{r}\cdot N_{t}\right)$
variables, which is huge and it will take a very long time for the
learning to converge. On the other hand, if we adopt the standard
Q-learning approach {[}27{]}-{[}29{]}, the domain of the Q-function
to be learned has $\left(1+S+N_{r}\cdot N_{t}+N_{t}\right)^{2}$ dimensions,
which is also prohibitively large. Secondly, the extended state space
$\left\{ \mathbf{S}_{k}\right\} $ has an infinite state space size.
Specifically, in standard Q-learning{[}27{]}-{[}29{]}, the size of
the state space is finite and a lookup table is used to store the
Q-value for the state-action pairs. However, for the extended state
$\mathbf{\mathbf{S}}_{k}=\left(\mathbf{x}_{k},\mathbf{H}_{k},\delta_{k}\right)$,
the wireless MIMO fading channel $\mathbf{H}_{k}$ switches in a continuous
state space with \textsl{uncountably many realizations}. As a result,
the size of the extended state space $\left\{ \mathbf{S}_{k}\right\} $
is \textsl{uncountably infinite}, which makes the associated learning
of the Q-function far more complicated.

As a result, instead of directly working on the standard \textsl{Bellman
optimality equation }(\ref{eq: standard Bellman Eq}), we derive a
\textsl{reduced-state Bellman equation} from (\ref{eq: standard Bellman Eq})
using the i.i.d. property of the channel state $\mathbf{H}_{k}$ and
the controller random access state $\delta_{k}$.

\begin{thm}
(\textsl{Reduced-State Bellman Optimality Equation)} \label{Thm: reduced-state-Bellman-Eq}
If there exists a pair of $\left(\widetilde{\theta},\widetilde{V}\left(\mathbf{x}_{k}\right)\right)$
that solves the following \textsl{reduced-state Bellman optimality
equation:}
\begin{align}
 & \widetilde{\theta}+\widetilde{V}\left(\mathbf{x}_{k}\right)=\nonumber \\
 & \mathbb{E}_{\left\{ \mathbf{H}_{k},\delta_{k}\right\} }\bigg[\min_{\mathbf{u}_{k}}\bigg[r\left(\mathbf{\mathbf{S}}_{k},\mathbf{u}_{k}\right)+\mathbb{E}\left[\left.V\left(\mathbf{S}_{k+1}\right)\right|\mathbf{\mathbf{S}}_{k},\mathbf{u}_{k}\right]\bigg]\bigg],\forall\mathbf{x}_{k},\label{eq: reduced-State-Bellman Eq}
\end{align}
then:
\end{thm}

\begin{itemize}
\item $\widetilde{\theta}=\theta=\min_{\mathbf{u}_{k}}\limsup_{K\rightarrow\infty}\frac{1}{K}\sum_{k=0}^{K}\mathbb{E}\left[r\left(\mathbf{\mathbf{S}}_{k},\mathbf{u}_{k}\right)\right]$
is the optimal average cost for Problem \ref{problem-1}.
\item $\widetilde{V}\left(\mathbf{x}_{k}\right)=\mathbb{E}\left[\left.V\left(\mathbf{S}_{k}\right)\right|\mathbf{x}_{k}\right]$
is the optimal \textsl{reduced state value function}.
\item The optimal control policy for Problem \ref{problem-1} is given by
$\mathbf{u}_{k}^{*}$, which attains the minimum of the R.H.S. of
(\ref{eq: reduced-State-Bellman Eq}).
\end{itemize}
\begin{IEEEproof}
Please see Appendix \ref{subsec:Proof-of-Thm1-Thm2}.
\end{IEEEproof}

Compared with the standard \textsl{Bellman optimality equation }(\ref{eq: standard Bellman Eq}),
the \textsl{reduced state Bellman optimality equation }(\ref{eq: reduced-State-Bellman Eq})\textsl{
}involves a\textsl{ reduced state value function} $\widetilde{V}\left(\mathbf{x}_{k}\right)$,
which is a function of the PSI $\mathbf{x}_{k}$ only. The number
of variables in the domain of $\widetilde{V}\left(\mathbf{x}_{k}\right)$
is reduced to $S$. As such, learning the reduced state value function
$\widetilde{V}\left(\mathbf{x}_{k}\right)$ is much easier than the
original value function $V\left(\mathbf{S}_{k}\right)$. In the next
sections, we shall exploit the specific structure of the \textsl{reduced
state Bellman optimality equation }(\ref{eq: reduced-State-Bellman Eq})
and characterize the existence and uniqueness of both the optimal
reduced state value function $\widetilde{V}\left(\mathbf{x}_{k}\right)$
and the optimal control action $\mathbf{u}_{k}^{*}$.

\subsection{Structural Properties of the Reduced State Value Function and Control
Action}

The structural properties of the solution to the \textsl{reduced state
Bellman optimality equation} (\ref{eq: reduced-State-Bellman Eq})
are important because they can provide potential opportunities to
simplify (\ref{eq: reduced-State-Bellman Eq}) and to develop the
online learning algorithm for reduced state value function $\widetilde{V}\left(\mathbf{x}_{k}\right)$
and the optimal control action $\mathbf{u}_{k}^{*}$. We summarize
the main result in the following theorem.

\begin{thm}
(\textsl{Structural Properties of the Solution to the Reduced-State
Bellman Optimality Equation)} \label{Thm: structure-property-reduced-state-Bellman}
If there exists a pair of $\left(\widetilde{\theta},\widetilde{V}\left(\mathbf{x}_{k}\right)\right)$
that solves the \textsl{reduced state Bellman optimality equation
}(\ref{eq: reduced-State-Bellman Eq}), then:
\end{thm}

\begin{itemize}
\item $\widetilde{V}\left(\mathbf{x}_{k}\right)$ is the optimal \textsl{reduced
state value function}. Moreover, $\widetilde{V}\left(\mathbf{x}_{k}\right)$
is quadratic w.r.t. the plant state $\mathbf{x}_{k}$ and is in the
form of 
\begin{align}
 & \widetilde{V}\left(\mathbf{x}_{k}\right)=\mathbf{x}_{k}^{T}\mathbf{P}\mathbf{x}_{k},\label{eq:opt value function}
\end{align}
where $\mathbf{P}\in\mathbb{S}_{+}^{S}$ is a constant positive definite
matrix.
\item $\widetilde{\theta}$ is the optimal average cost and is in the form
of $\widetilde{\theta}=\mathrm{Tr}\left(\mathbf{M}+\mathbf{P}\mathbf{W}+\mathbf{B}^{T}\mathbf{P}\mathbf{B}\right)$.
\item The optimal control action $\mathbf{u}_{k}^{*}=\left(\mathrm{\Omega}^{k}\right)^{*}\left(\mathbf{\mathbf{S}}_{k}\right)$
for a given state realization $\mathbf{\mathbf{S}}_{k}$ is a linear
state feedback control law
\begin{align}
 & \mathbf{u}_{k}^{*}=\left(\mathrm{\Omega}^{k}\right)^{*}\left(\mathbf{\mathbf{S}}_{k}\right)=-(\delta_{k}\mathbf{H}_{k}^{T}\mathbf{B}^{T}\mathbf{P}\mathbf{B}\mathbf{H}_{k}\nonumber \\
 & +\delta_{k}\mathbf{H}_{k}^{T}\mathbf{M}\mathbf{H}_{k}+\mathbf{R})^{-1}\mathbf{H}_{k}^{T}\mathbf{B}^{T}\mathbf{P}\mathbf{A}\mathbf{\mathbf{x}}_{k}.\label{eq:optimal control action}
\end{align}
\end{itemize}
\begin{IEEEproof}
Please see Appendix C.
\end{IEEEproof}

Based on the structural properties of the \textsl{reduced state value
function} $\widetilde{V}\left(\mathbf{x}_{k}\right)$, the optimal
average cost $\widetilde{\theta}$, and the optimal control policy
$\left(\mathrm{\Omega}^{k}\right)^{*}$ (\ref{eq:optimal control action})
in Theorem \ref{Thm: structure-property-reduced-state-Bellman}, it
follows that the existence and uniqueness of $\widetilde{V}\left(\mathbf{x}_{k}\right)$,
$\widetilde{\theta}$ and $\left(\mathrm{\Omega}^{k}\right)^{*}$
is equivalent to the existence and uniqueness of the kernel $\mathbf{P}$
of the \textsl{reduced state value function }$\widetilde{V}\left(\mathbf{x}_{k}\right)$.
Moreover, substituting the structural forms of $\widetilde{V}\left(\mathbf{x}_{k}\right)$,
$\widetilde{\theta}$ and $\left(\mathrm{\Omega}^{k}\right)^{*}$
in Theorem \ref{Thm: structure-property-reduced-state-Bellman} into
(\ref{eq: reduced-State-Bellman Eq}), it follows that the \textsl{reduced
state Bellman optimality equation} (\ref{eq: reduced-State-Bellman Eq})
can be further simplified as:
\begin{align}
 & \mathbf{\mathbf{x}}_{k}^{T}\mathbf{P}\mathbf{\mathbf{x}}_{k}=\mathbf{\mathbf{x}}_{k}^{T}(\mathbf{A}^{T}\mathbf{P}\mathbf{A}-\mathbb{E}[\delta_{k}\mathbf{A}^{T}\mathbf{P}\mathbf{B}\mathbf{H}_{k}(\delta_{k}\mathbf{H}_{k}^{T}\mathbf{B}^{T}\mathbf{P}\mathbf{B}\mathbf{H}_{k}\nonumber \\
 & +\delta_{k}\mathbf{H}_{k}^{T}\mathbf{M}\mathbf{H}_{k}+\mathbf{R})^{-1}\mathbf{H}_{k}^{T}\mathbf{B}^{T}\mathbf{P}\mathbf{A}]+\mathbf{Q})\mathbf{\mathbf{x}}_{k},\forall\mathbf{x}_{k}.\label{eq:equivalent reduced state bellman eq}
\end{align}

Since equation (\ref{eq:equivalent reduced state bellman eq}) must
be satisfied for all $\mathbf{x}_{k}$, we have the following lemma
on the existence and the uniqueness of $\widetilde{V}\left(\mathbf{x}_{k}\right)$,
$\theta$ and $\mathbf{u}_{k}^{*}$.

\begin{lemma}
\textsl{(Existence and Uniqueness of $\widetilde{V}\left(\mathbf{x}_{k}\right)$,
$\widetilde{\theta}$ and $\mathbf{u}_{k}^{*}$)}\label{lemma-existence-uniqueness}
If there exists a unique $\mathbf{P}\in\mathbb{S}_{+}^{S}$ such that
the following nonlinear matrix equation (NME) is satisfied:
\begin{align}
\mathbf{P}= & \mathbf{A}^{T}\mathbf{P}\mathbf{A}-\mathbb{E}[\delta_{k}\mathbf{A}^{T}\mathbf{P}\mathbf{B}\mathbf{H}_{k}(\delta_{k}\mathbf{H}_{k}^{T}\mathbf{B}^{T}\mathbf{P}\mathbf{B}\mathbf{H}_{k}\nonumber \\
 & +\delta_{k}\mathbf{H}_{k}^{T}\mathbf{M}\mathbf{H}_{k}+\mathbf{R})^{-1}\mathbf{H}_{k}^{T}\mathbf{B}^{T}\mathbf{P}\mathbf{A}]+\mathbf{Q},\label{eq: NME}
\end{align}
then the solution pair of $\left(\widetilde{\theta},\widetilde{V}\left(\mathbf{x}_{k}\right)\right)$
that solves the \textsl{reduced state Bellman optimality equation
}(\ref{eq: reduced-State-Bellman Eq}) and the associated optimal
control action $\mathbf{u}_{k}^{*}$ exists and is unique.
\end{lemma}

\begin{IEEEproof}
Please see Appendix C.
\end{IEEEproof}

As a result, we shall focus on the NME (\ref{eq: NME}), and analyze
the existence and uniqueness of $\mathbf{P}$ that satisfies (\ref{eq: NME}).

\subsection{Sufficient Condition for Existence and Uniqueness of Optimal Control
Action}

Note that in a special case of LTI systems, where the communication
channel between the remote controller and the actuator is static with
$\mathbf{H}_{k}=\mathbf{I}$ and $\delta_{k}=1$, the NME (\ref{eq: NME})
is reduced to the following standard \textsl{discrete-time algebraic
Riccati equation }(DARE):
\begin{align}
\mathbf{P}= & \mathbf{A}^{T}\mathbf{P}\mathbf{A}-\mathbf{A}^{T}\mathbf{P}\mathbf{B}(\mathbf{B}^{T}\mathbf{P}\mathbf{B}\nonumber \\
 & +\mathbf{M}+\mathbf{R})^{-1}\mathbf{B}^{T}\mathbf{P}\mathbf{A}+\mathbf{Q}.\label{eq: DARE}
\end{align}
It is well studied in existing literature {[}30{]} and {[}31{]}, that
if the LTI system is controllable, i.e., the pair $\left(\mathbf{A},\mathbf{B}\right)$
is controllable, then there exist a unique $\mathbf{P}$ that satisfies
the DARE (\ref{eq: DARE}). Moreover, the associated optimal control
action that solves the \textsl{reduced state Bellman optimality equation}
(\ref{eq: reduced-State-Bellman Eq}) is the certainty equivalent
controller $\mathbf{u}_{k}^{*}=-(\mathbf{B}^{T}\mathbf{P}\mathbf{B}+\mathbf{M}+\mathbf{R})^{-1}\mathbf{B}^{T}\mathbf{P}\mathbf{A}\mathbf{x}_{k}$
{[}31{]}.

However, in the presence of random access of the controller $\delta_{k}$
and random MIMO fading channel $\mathbf{H}_{k}$, the controllability
of the closed-loop control system may not be preserved even if the
pair $\left(\mathbf{A},\mathbf{B}\right)$ is controllable. As a result,
it is important to characterize the impacts of general random MIMO
fading channels and random access of the controller on the controllability
of the closed-loop system. This is formally summarized in the following
Lemma.

\begin{lemma}
\textsl{(Impacts of MIMO Fading Channel and Random Access on Closed-loop
Controllability)}\label{lemma-impacts-MIMO Fading} Let the the singular
value decomposition (SVD) of $\mathbf{B}\mathbf{B}^{T}$ be $\mathbf{B}\mathbf{B}^{T}=\mathbf{U}^{T}\boldsymbol{\boldsymbol{\Xi}}\mathbf{U}$
with the diagonal elements of $\boldsymbol{\boldsymbol{\Xi}}$ in
descending order, where $\mathbf{U}\in\mathbb{R}^{S\times S}$ is
an unitary matrix. Denote $\mathrm{rank}\left(\mathbf{B}\right)=\eta_{B}$.
Let the similarity transformation of $\mathbf{A}$ w.r.t. $\mathbf{U}$
be $\widetilde{\mathbf{A}}=\mathbf{U}^{T}\mathbf{A}\mathbf{U}$. Denote
the block-wise representation of $\widetilde{\mathbf{A}}$ as $\widetilde{\mathbf{A}}=\left[\begin{array}{cc}
\widetilde{\mathbf{A}}_{11} & \widetilde{\mathbf{A}}_{12}\\
\widetilde{\mathbf{A}}_{21} & \widetilde{\mathbf{A}}_{22}
\end{array}\right]$, where $\widetilde{\mathbf{A}}_{11}\in\mathbb{R}^{\eta_{B}\times\eta_{B}}$
is the $\eta_{B}$-th order leading principal submatrix of matrix
$\widetilde{\mathbf{A}}$, $\widetilde{\mathbf{A}}_{12}\in\mathbb{R}^{\eta_{B}\times\left(S-\eta_{B}\right)}$,
$\widetilde{\mathbf{A}}_{12}\in\mathbb{R}^{\left(S-\eta_{B}\right)\times\eta_{B}}$,
and $\widetilde{\mathbf{A}}_{22}\in\mathbb{R}^{\left(S-\eta_{B}\right)\times\left(S-\eta_{B}\right)}$
are constant matrices. Assume the pair $\left(\mathbf{A},\mathbf{B}\right)$
is controllable, the impacts of $\mathbf{H}_{k}$ and $\delta_{k}$
on the controllability of $\left(\mathbf{A},\delta_{k}\mathbf{B}\mathbf{H}_{k}\right)$
are given by:
\end{lemma}

\begin{itemize}
\item \textsl{(a) Almost Sure Controllability}: If one of the following
three conditions \textsl{(a.1)}, \textsl{(a.2)} and \textsl{(a.3)}
is satisfied, then the closed-loop system is almost surely controllable,
i.e., $\left(\mathbf{A},\delta_{k}\mathbf{B}\mathbf{H}_{k}\right)$
is controllable w.p.1. for any timeslot $k$.
\begin{itemize}
\item \textsl{(a.1)} $N_{t}\geq S$ and $\mathrm{Pr}\left(\delta_{k}=1\right)=1,\forall k\geq0$$;$
\item \textsl{(a.2)} $\eta_{B}\leq N_{t}<S$, $\mathrm{Pr}\left(\delta_{k}=1\right)=1,\forall k\geq0$,
and the pair $\left(\widetilde{\mathbf{A}}_{22},\widetilde{\mathbf{A}}_{12}^{T}\right)$
is controllable.
\item \textsl{(a.3)} $N_{t}<\eta_{B}$, $\mathrm{Pr}\left(\delta_{k}=1\right)=1,\forall k\geq0$,
and $\mathrm{Rank}\left(\widetilde{\mathbf{A}}-\lambda\mathbf{I}\right)>\left(S-\eta_{B}+N_{t}\right)$,
$\forall\lambda\in\mathbb{C}$.
\end{itemize}
\item \textsl{(b)} \textsl{Intermittent Controllability}: If one of the
following three conditions \textsl{(b.1)}, \textsl{(b.2)} and \textsl{(b.3)}
is satisfied, then the closed-loop system is intermittent controllable,
i.e., $\left(\mathbf{A},\delta_{k}\mathbf{B}\mathbf{H}_{k}\right)$
is almost surely controllable at the timeslot $k$ when $\delta_{k}=1$,
and $\left(\mathbf{A},\delta_{k}\mathbf{B}\mathbf{H}_{k}\right)$
is uncontrollable at the timeslot $k$ when $\delta_{k}=0$.
\begin{itemize}
\item \textsl{(b.1)} $N_{t}\geq S$ and $0<\mathrm{Pr}\left(\delta_{k}=1\right)<1$
, $\forall k\geq0$;
\item \textsl{(b.2)} $\eta_{B}\leq N_{t}<S$, $0<\mathrm{Pr}\left(\delta_{k}=1\right)<1,\forall k\geq0$,
and the pair $\left(\widetilde{\mathbf{A}}_{22},\widetilde{\mathbf{A}}_{12}^{T}\right)$
is controllable.
\item \textsl{(b.3)} $N_{t}<\eta_{B}$, $0<\mathrm{Pr}\left(\delta_{k}=1\right)<1,\forall k\geq0$,
and $\mathrm{Rank}\left(\widetilde{\mathbf{A}}-\lambda\mathbf{I}\right)>\left(S-\eta_{B}+N_{t}\right)$,
$\forall\lambda\in\mathbb{C}$.
\end{itemize}
\item \textsl{(c) Almost Sure Uncontrollability}: If either of the following
two conditions \textsl{(c.1)} or \textsl{(c.2)} is satisfied, then
the closed-loop system is almost surely uncontrollable, i.e., $\left(\mathbf{A},\delta_{k}\mathbf{B}\mathbf{H}_{k}\right)$
is uncontrollable w.p.1. at any timeslot $k$.
\begin{itemize}
\item \textsl{(c.1)} $\eta_{B}\leq N_{t}<S$ and the pair $\left(\widetilde{\mathbf{A}}_{22},\widetilde{\mathbf{A}}_{12}^{T}\right)$
is uncontrollable.
\item \textsl{(c.2)} $N_{t}<\eta_{B}$ and $\exists\lambda\in\mathbb{C}$
such that $\mathrm{Rank}\left(\widetilde{\mathbf{A}}-\lambda\mathbf{I}\right)\leq\left(S-\eta_{B}+N_{t}\right).$
\end{itemize}
\end{itemize}
\begin{IEEEproof}
Please see Appendix D.
\end{IEEEproof}

The boundaries that distinguish almost sure controllability, intermittent
controllability and almost sure uncontrollability are visualized in
the following Fig. \ref{fig: controllability-boundaries}.

\begin{figure}[tbh]
\begin{centering}
\includegraphics[clip,width=1\columnwidth]{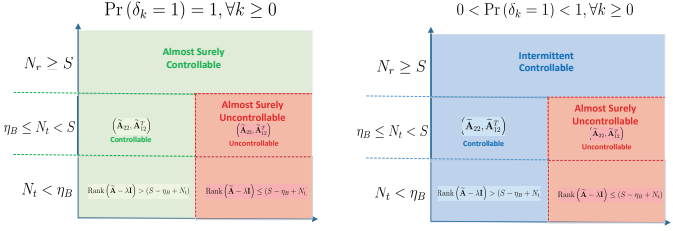}
\par\end{centering}
\caption{\label{fig: controllability-boundaries}Illustration of the regions
of almost sure controllability, intermittent controllability and almost
sure uncontrollability.}
\end{figure}

Due to the intermittent controllability and almost sure uncontrollability
caused by the random access of the controller $\delta_{k}$ and random
MIMO fading channel $\mathbf{H}_{k}$, the existing approaches for
analyzing the solution of the DARE (\ref{eq: DARE}) are not applicable
to the NME (\ref{eq: NME}). Moreover, the expectation w.r.t. both
$\delta_{k}$ and $\mathbf{H}_{k}$ on the R.H.S. of the NME (\ref{eq: NME})
does not have closed-form expression due to the high nonlinearity
of the matrix inversion. We shall address this challenge via exploiting
the underlying structure of the NME (\ref{eq: NME}).

We observe that if kernel $\mathbf{P}$ of the \textsl{reduced state
value function }$\widetilde{V}\left(\mathbf{x}_{k}\right)$ is such
that $\delta_{k}\mathbf{H}_{k}^{T}\mathbf{B}^{T}\mathbf{P}=\boldsymbol{0}$,
then the NME (\ref{eq: NME}) can be substantially simplified as $\mathbf{P}=\mathbf{A}\mathbf{P}\mathbf{A}^{T}+\mathbf{W}.$
The condition for the existence and uniqueness for such $\mathbf{P}$
is readily established in closed-form as $\rho\left(\mathbf{A}\right)<1$
{[}32{]}. Motivated by this fact, given $\delta_{k}$ and $\mathbf{H}_{k}$,
for any arbitrary $\mathbf{P}$, we aim at decomposing $\mathbf{P}$
into a sum of two PSD components $\boldsymbol{\mathbf{P}}_{k}^{c}$
and $\boldsymbol{\mathbf{P}}_{k}^{uc}$ such that $\delta_{k}\mathbf{H}_{k}^{T}\mathbf{B}^{T}\boldsymbol{\mathbf{P}}_{k}^{uc}=\boldsymbol{0}$.
This leads to the following definition of the controllable and uncontrollable
cones of PSD matrices.

\begin{definitn}
\textsl{(Controllable and Uncontrollable PSD Cones)} Given a certain
realization of $\left(\delta_{k},\mathbf{H}_{k}\right)$, the controllable
PSD cone $\mathcal{C}^{c}$ and uncontrollable PSD cone $\mathcal{C}^{uc}$
associated with $\delta_{k}\mathbf{B}\mathbf{H}_{k}$ is defined by
\begin{align}
 & \mathcal{C}^{c}=\left\{ \left.\mathbf{T}\in\mathbb{S}^{S}\right|\mathrm{\text{ker}}\left(\delta_{k}\mathbf{B}\mathbf{H}_{k}\mathbf{H}_{k}^{T}\mathbf{B}^{T}\mathbf{T}\right)=\mathrm{\text{ker}}\left(\mathbf{T}\right)\right\} ;\\
 & \mathcal{C}^{uc}=\left\{ \left.\mathbf{T}\in\mathbb{S}^{S}\right|\delta_{k}\mathbf{B}\mathbf{H}_{k}\mathbf{H}_{k}^{T}\mathbf{B}^{T}\mathbf{T}=\mathbf{0}\right\} .
\end{align}
\end{definitn}

Let the SVD of $\delta_{k}\mathbf{B}\mathbf{H}_{k}\left(\delta_{k}\mathbf{H}_{k}^{T}\mathbf{M}\mathbf{H}_{k}+\mathbf{R}\right)^{-1}\mathbf{H}_{k}^{T}\mathbf{B}^{T}$
be
\begin{align}
 & \delta_{k}\mathbf{B}\mathbf{H}_{k}\left(\delta_{k}\mathbf{H}_{k}^{T}\mathbf{M}\mathbf{H}_{k}+\mathbf{R}\right)^{-1}\mathbf{H}_{k}^{T}\mathbf{B}^{T}=\mathbf{V}_{k}^{T}\boldsymbol{\Lambda}_{k}\mathbf{V}_{k},
\end{align}
with the diagonal elements of $\boldsymbol{\Lambda}_{k}$ in descending
order. Let $\mathrm{rank}\left(\delta_{k}\mathbf{B}\mathbf{H}_{k}\mathbf{H}_{k}^{T}\mathbf{B}^{T}\right)=\gamma_{k}$
and $\boldsymbol{\Pi}_{k}=\left[\begin{array}{cc}
\mathbf{I}_{\gamma_{k}} & \mathbf{0}\\
\mathbf{0} & \mathbf{0}
\end{array}\right]_{S\times S}.$ The closed-form controllable and uncontrollable PSD cone decomposition
of $\mathbf{P}$ is characterized by the following Theorem.

\begin{thm}
\textsl{(Closed-form PSD Cone Decomposition of} $\mathbf{P}$\textsl{)}
\label{Thm-PSD decomposition}Given a certain realization of $\left(\delta_{k},\mathbf{H}_{k}\right)$,
the kernel of the value function $\mathbf{P}$ can be decomposed into
two parts as $\mathbf{P}=\boldsymbol{\mathbf{P}}_{k}^{c}+\boldsymbol{\mathbf{P}}_{k}^{uc}$,
where $\boldsymbol{\mathbf{P}}_{k}^{c}\in\mathcal{C}^{c}$ and $\boldsymbol{\mathbf{P}}_{k}^{uc}\in\mathcal{C}^{uc}$.
The closed-form expressions of $\boldsymbol{\mathbf{P}}_{k}^{c}$
and $\boldsymbol{\mathbf{P}}_{k}^{uc}$ are given by:
\begin{align}
\boldsymbol{\mathbf{P}}_{k}^{c}= & \mathbf{V}_{k}^{T}\left[\begin{array}{cc}
\left(\mathbf{V}_{k}\mathbf{P}\mathbf{V}_{k}^{T}\right)_{\gamma_{k}} & \left(\mathbf{V}_{k}\mathbf{P}\mathbf{V}_{k}^{T}\right)_{\gamma_{k}}\boldsymbol{\Sigma}_{k}\\
\boldsymbol{\Sigma}_{k}^{T}\left(\mathbf{V}_{k}\mathbf{P}\mathbf{V}_{k}^{T}\right)_{\gamma_{k}} & \boldsymbol{\Sigma}_{k}^{T}\left(\mathbf{V}_{k}\mathbf{P}\mathbf{V}_{k}^{T}\right)_{\gamma_{k}}\boldsymbol{\Sigma}_{k}
\end{array}\right]\mathbf{V}_{k};\label{eq:Pk-c}\\
\boldsymbol{\mathbf{P}}_{k}^{uc}= & \mathbf{V}_{k}^{T}\left(\mathbf{I}-\boldsymbol{\Pi}_{k}\right)\mathbf{V}_{k}\mathbf{P}\mathbf{V}_{k}^{T}\left(\mathbf{I}-\boldsymbol{\Pi}_{k}\right)\mathbf{V}_{k}\nonumber \\
 & -\mathbf{V}_{k}^{T}\mathrm{diag}\left(\mathbf{0}_{\gamma_{k}},\boldsymbol{\Sigma}_{k}^{T}\left(\mathbf{V}_{k}\mathbf{P}\mathbf{V}_{k}^{H}\right)_{\gamma_{k}}\boldsymbol{\Sigma}_{k}\right)\mathbf{V}_{k},\label{eq:Pk-uc}
\end{align}
where $\boldsymbol{\Sigma}_{k}=\left(\mathbf{V}_{k}\mathbf{P}\mathbf{V}_{k}^{T}\right)_{\gamma_{k}}^{-1}\left(\mathbf{V}_{k}\mathbf{P}\mathbf{V}_{k}^{T}\right)_{\left(1:\gamma_{k};\gamma_{k}+1:S\right)}$.
\end{thm}

\begin{IEEEproof}
Please see Appendix E.
\end{IEEEproof}

Utilizing the controllable and uncontrollable PSD cone decomposition
of $\mathbf{P}$, the NME (\ref{eq: NME}) can also be decomposed
into two parts with each part depending solely on $\boldsymbol{\mathbf{P}}_{k}^{c}$
and $\boldsymbol{\mathbf{P}}_{k}^{uc}$. This is formally stated in
the following Lemma.

\begin{lemma}
\textsl{(Decomposition of the NME of }$\boldsymbol{\mathbf{P}}$\textsl{)}
\label{lemma-decomposed-NME}Applying the controllable and uncontrollable
PSD cone decomposition of $\mathbf{P}$ in Theorem \ref{Thm-PSD decomposition},
the NME of $\mathbf{P}$ in (\ref{eq: NME}) can be represented in
a more fine-grained form given by
\begin{align}
 & \mathbf{P}=\mathbf{Q}+\underbrace{\mathbf{A}^{T}\mathbb{E}\left[\boldsymbol{\mathbf{P}}_{k}^{uc}\right]\mathbf{A}}_{\emph{depends solely on \ensuremath{\boldsymbol{\mathbf{P}}_{k}^{uc}}}}+\nonumber \\
 & \underbrace{\mathbf{A}^{T}\mathbb{E}[\boldsymbol{\mathbf{P}}_{k}^{c}\mathbf{A}-\delta_{k}\mathbf{A}^{T}\boldsymbol{\mathbf{P}}_{k}^{c}\mathbf{B}\mathbf{H}_{k}(\delta_{k}\mathbf{H}_{k}^{T}\mathbf{B}^{T}\boldsymbol{\mathbf{P}}_{k}^{c}\mathbf{B}\mathbf{H}_{k}}_{\emph{depends solely on \ensuremath{\boldsymbol{\mathbf{P}}_{k}^{c}}}}\nonumber \\
 & \underbrace{+\delta_{k}\mathbf{H}_{k}^{T}\mathbf{M}\mathbf{H}_{k}+\mathbf{R})^{-1}\mathbf{H}_{k}^{T}\mathbf{B}^{T}\boldsymbol{\mathbf{P}}_{k}^{c}]\mathbf{A}}_{\emph{depends solely on \ensuremath{\boldsymbol{\mathbf{P}}_{k}^{c}}}}.\label{eq:fine-grained-NME}
\end{align}
\end{lemma}

\begin{IEEEproof}
Please see Appendix E.
\end{IEEEproof}

Compared with the original NME of $\mathbf{P}$ in (\ref{eq: NME})
, the fine-grained decomposed NME (\ref{eq:fine-grained-NME}) is
more informative and reveals key insights into the existence of $\mathbf{P}$.
Intuitively, when $\boldsymbol{\mathbf{P}}_{k}^{c}$ is dominant (i.e.,
$\boldsymbol{\mathbf{P}}_{k}^{c}\gg\boldsymbol{\mathbf{P}}_{k}^{uc}$),
the existence of $\boldsymbol{\mathbf{P}}$ can be obtained via the
monotonicity and boundedness of the R.H.S. of (\ref{eq:fine-grained-NME}).
However, when $\boldsymbol{\mathbf{P}}_{k}^{uc}$ is dominant (i.e.,
$\boldsymbol{\mathbf{P}}_{k}^{uc}\gg\boldsymbol{\mathbf{P}}_{k}^{c}$),
$\boldsymbol{\mathbf{P}}$ tends to be unstable because the R.H.S.
of (\ref{eq:fine-grained-NME}) is dominated by $\mathbf{A}^{T}\mathbb{E}\left[\boldsymbol{\mathbf{P}}_{k}^{uc}\right]\mathbf{A}$,
and the $\boldsymbol{\mathbf{P}}$ that satisfies (\ref{eq:fine-grained-NME})
may not exist if the dynamic plant is unstable, i.e., $\rho\left(\mathbf{A}\right)>1$.
Therefore, $\boldsymbol{\mathbf{P}}_{k}^{c}$ is favorable, whereas
$\boldsymbol{\mathbf{P}}_{k}^{uc}$ is unfavorable for the existence
of $\boldsymbol{\mathbf{P}}$ that satisfies (\ref{eq:fine-grained-NME}).
As a result, instead of analyzing the original NME (\ref{eq: NME})
in a brute-force manner, we shall focus on the analysis of the existence
and uniqueness of $\mathbf{P}$ based on the fine-grained decomposed
NME (\ref{eq:fine-grained-NME}). The closed-form sufficient condition
for the existence and uniqueness of $\mathbf{P}$ that satisfies (\ref{eq:fine-grained-NME})
is summarized in the following Theorem \ref{Thm: suff-condition-existence-uniqueness}.

\begin{thm}
(\textsl{Sufficient Condition for the Existence and the Uniqueness
of Optimal Control)} \label{Thm: suff-condition-existence-uniqueness}If
one of the three conditions \textsl{(a.1)}, \textsl{(a.2)} and \textsl{(a.3)}
in Lemma \ref{lemma-impacts-MIMO Fading} is satisfied, or the following
condition (\ref{eq:suff-condition-existence}): 
\begin{align}
 & \left\Vert \mathbb{E}\left[\mathbf{A}^{T}\mathbf{V}_{k}^{T}\left(\mathbf{I}-\boldsymbol{\Pi}_{k}\right)\mathbf{V}_{k}\mathbf{A}\right]\right\Vert <1,\label{eq:suff-condition-existence}
\end{align}
is satisfied, then the solution $\mathbf{P}$ to the NME (\ref{eq: NME})
exists and is unique. Moreover, the optimal control action $\mathbf{u}^{*}\left(\mathbf{x}_{k}\right)$
that solves the infinite horizon ergodic control Problem \ref{problem-1}
exists and is unique, and is given by (\ref{eq:optimal control action}).
\end{thm}

\begin{IEEEproof}
Please see Appendix F.
\end{IEEEproof}

Theorem \ref{Thm: suff-condition-existence-uniqueness} delivers a
key message that provided the closed-form sufficient condition (\ref{eq:suff-condition-existence})
is satisfied, the optimal control action still exists and is unique
even if the closed-loop control system suffers from intermittent controllability
or almost sure uncontrollability.

\section{Simultaneous Learning of the Value Function and Control Over MIMO
Fading Channels}

In this section, we shall propose an online and autonomous learning
algorithm that can simultaneously learn both the optimal control action
$\mathbf{u}^{*}\left(\mathbf{S}_{k}\right)$ and the associated optimal
value function $\widetilde{V}\left(\mathbf{x}_{k}\right)$ on the
fly based on the state observations $\mathbf{x}_{k}$ only. The proposed
online solution is implemented at the remote controller and has fast
convergence.

\subsection{Simultaneous Learning of Value Function and Optimal Control Action}

Note that the NME (\ref{eq: NME}) is an algebraic equation with unknown
variable $\mathbf{P}$. Thus, we shall utilize the stochastic approximation
theory to construct an online learning algorithm to estimate the unknown
$\mathbf{P}$ based on the algebraic equation (\ref{eq: NME}). The
estimated variable $\mathbf{P}$ can then be applied to obtain the
optimal value function $\widetilde{V}\left(\mathbf{x}_{k}\right)$
and the optimal control action $\mathbf{u}^{*}\left(\mathbf{S}_{k}\right)$
simultaneously using (\ref{eq:opt value function}) and (\ref{eq:optimal control action}),
respectively, in Theorem \ref{Thm: structure-property-reduced-state-Bellman}.

We first rewrite the NME (\ref{eq: NME}) into the standard form of
$f(\mathbf{P})=\mathbf{0}$ and apply the stochastic approximation
(SA) technique to estimate the root of the equation. Specifically,
\begin{align}
 & f(\mathbf{P})=\mathbf{A}^{T}\mathbf{P}\mathbf{A}-\mathbb{E}[\delta_{k}\mathbf{A}^{T}\mathbf{P}\mathbf{B}\mathbf{H}_{k}(\delta_{k}\mathbf{H}_{k}^{T}\mathbf{B}^{T}\mathbf{P}\mathbf{B}\mathbf{H}_{k}\nonumber \\
 & +\delta_{k}\mathbf{H}_{k}^{T}\mathbf{M}\mathbf{H}_{k}+\mathbf{R})^{-1}\mathbf{H}_{k}^{T}\mathbf{B}^{T}\mathbf{P}\mathbf{A}]-\mathbf{P}+\mathbf{Q}.\label{eq: SA-f(P)}
\end{align}

\begin{itemize}
\item \textbf{Online Learning of $\mathbf{P}^{*}$: }To obtain the root
of $f(\mathbf{P})=\mathbf{0}$, we can apply the stochastic approximation
iteration,
\begin{align}
\mathbf{P}_{k+1} & =\mathbf{P}_{k}+\alpha_{k}\widehat{f}\left(\mathbf{P}_{k}\right),\forall k\geq0,\label{eq: prps-SA}
\end{align}
where $\mathbf{P}_{0}\in\mathbb{S}_{+}^{S}$ is a bounded constant
positive definite matrix, $\left\{ \alpha_{k},k\geq0\right\} $ is
the step-size sequence satisfying
\begin{align}
\sum_{k=1}^{\infty}\alpha_{k}=\infty, & \sum_{k=1}^{\infty}\alpha_{k}^{2}<\infty,
\end{align}
and $\widehat{f}\left(\mathbf{P}_{k}\right)$ is an unbiased estimator
of $f(\mathbf{P}_{k})$, i.e., $f\left(\mathbf{P}_{k}\right)=\mathbb{E}\left[\left.\widehat{f}\left(\mathbf{P}_{k}\right)\right|\mathbf{P}_{k}\right]$,
and is given by
\begin{align}
\widehat{f}\left(\mathbf{P}_{k}\right)= & \mathbf{A}^{T}\mathbf{P}_{k}\mathbf{A}-\delta_{k}\mathbf{A}^{T}\mathbf{P}_{k}\mathbf{B}\mathbf{H}_{k}(\delta_{k}\mathbf{H}_{k}^{T}\mathbf{B}^{T}\mathbf{P}_{k}\mathbf{B}\mathbf{H}_{k}\nonumber \\
 & +\delta_{k}\mathbf{H}_{k}^{T}\mathbf{M}\mathbf{H}_{k}+\mathbf{R})^{-1}\mathbf{H}_{k}^{T}\mathbf{B}^{T}\mathbf{P}_{k}\mathbf{A}-\mathbf{P}_{k}+\mathbf{Q}.\label{eq: unbiased estimator f-hat}
\end{align}
\item \textbf{Online learning of Control Action $\mathbf{u}_{k}$:} At the
$k$-th timeslot, the reduced state value function and the control
action can be computed based on $\mathbf{P}_{k}$:
\begin{align}
\widetilde{V}_{k}\left(\mathbf{x}_{k}\right)= & \mathbf{x}_{k}^{T}\mathbf{P}_{k}\mathbf{x}_{k},\label{eq: prps-Vk-iter}\\
\mathbf{u}_{k}= & -(\delta_{k}\mathbf{H}_{k}^{T}\mathbf{B}^{T}\mathbf{P}_{k}\mathbf{B}\mathbf{H}_{k}+\delta_{k}\mathbf{H}_{k}^{T}\mathbf{M}\mathbf{H}_{k}\nonumber \\
 & +\mathbf{R})^{-1}\mathbf{H}_{k}^{T}\mathbf{B}^{T}\mathbf{P}_{k}\mathbf{A}\mathbf{\mathbf{x}}_{k}.\label{eq: prps-u-iter}
\end{align}
\end{itemize}
In the iteration (\ref{eq: prps-SA}), only the realizations of $\mathbf{H}_{k}$
and $\delta_{k}$ will be required. The random access state $\delta_{k}$
is locally available at the remote controller. The wireless fading
channel realization $\mathbf{H}_{k}$ can be obtained by standard
channel estimation at the actuator based on the received pilot symbols
from the remote controller and channel feedback to the controller\footnote{In the LTE standard {[}34{]}, besides transmitting the control action
$\mathbf{u}_{k}$, the remote controller also transmits a pilot symbol
$\mathbf{L}\in\mathbb{R}^{N_{t}\times L}$ ($L>N_{t}$) on the PUSCH
data frame to the actuator at each timeslot. The actuator obtains
the MIMO fading channel realization $\mathbf{H}_{k}$ based the received
pilot signal $\mathbf{Y}_{k}^{L}=\mathbf{H}_{k}\mathbf{L}+\mathbf{V}_{k}^{L}$,
where $\mathbf{V}_{k}^{L}\in\mathbb{R}^{N_{r}\times L}$ is the additive
channel noise. The actuator then feeds back the CSI $\mathbf{H}_{k}$
to the controller.}.

The following lemma summarizes several key properties of the proposed
stochastic approximation iteration (\ref{eq: prps-SA}).

\begin{lemma}
(\textsl{Properties of the SA Iteration (\ref{eq: prps-SA}))} \label{Lemma: properties-iter}
\end{lemma}

\begin{itemize}
\item \textsl{(Lipschitz Continuity)} The matrix-valued function $f(\mathbf{P})$
is Lipschitz continuous with Lipschitz constant $\left\Vert \mathbf{A}\right\Vert ^{2}$,
i.e., $\left\Vert f(\mathbf{P}^{(1)})-f(\mathbf{P}^{(2)})\right\Vert \leq\left(1+\left\Vert \mathbf{A}\right\Vert ^{2}\right)\left\Vert \mathbf{P}^{(1)}-\mathbf{P}^{(2)}\right\Vert ,$
$\forall\mathbf{P}^{(1)},\mathbf{P}^{(2)}\in\mathbb{S}_{+}^{S}$.
\item \textsl{(Martingale Difference Noise)} Denote the estimation noise
of $f(\mathbf{P}_{k})$ in (\ref{eq: prps-SA}) as $\mathbf{N}_{k}=\left(\widehat{f}\left(\mathbf{P}_{k}\right)-f\left(\mathbf{P}_{k}\right)\right)$.
The sequence $\left\{ \mathbf{N}_{k},k\geq0\right\} $ is a martingale
difference sequence w.r.t. the filtration $\left\{ \mathcal{F}_{k}\triangleq\sigma\left(\mathbf{P}_{0},\delta_{1},\mathbf{H}_{1},\ldots,\delta_{k},\mathbf{H}_{k}\right)\right\} $,
i.e., $\mathbb{E}\left[\left.\mathbf{N}_{k+1}\right|\mathcal{F}_{k}\right]=\mathbf{0}_{S\times S},\forall k>0.$
\item \textsl{(Square Integrability)} $\left\{ \mathbf{N}_{k},k\geq0\right\} $
are square-integrable with $\mathbb{E}\left[\left.\left\Vert \mathbf{N}_{k+1}\right\Vert ^{2}\right|\mathcal{F}_{k}\right]\leq2\left\Vert \mathbf{A}\right\Vert ^{2}\left(1+\left\Vert \mathbf{P}_{k}\right\Vert ^{2}\right),\forall k>0.$
\end{itemize}
\begin{IEEEproof}
Please see Appendix G.
\end{IEEEproof}

In the next subsection, we shall focus on the convergence analysis
of the stochastic approximation iteration (\ref{eq: prps-SA}).

\subsection{Virtual Fixed-Point Process}

The ordinary differential equation (ODE) approach {[}35{]} serves
as a powerful tool for analyzing the limiting behaviors of the stochastic
approximation iteration (\ref{eq: prps-SA}). Specifically, rearranging
(\ref{eq: prps-SA}), it follows that
\begin{align}
\frac{\mathbf{P}_{k+1}-\mathbf{P}_{k}}{\alpha_{k}} & =f\left(\mathbf{P}_{k}\right)+\mathbf{N}_{k},\forall k\geq0.\label{eq:Rearranged-SA}
\end{align}
Intuitively, when $\alpha_{k}$ is sufficiently small, the nonlinear
difference equation (\ref{eq:Rearranged-SA}) can be approximated
by the following ODE:
\begin{align}
\dot{\mathbf{P}}\left(t\right) & =f\left(\mathbf{P}\left(t\right)\right),\mathbf{P}\left(t\right)=\mathbf{P}_{0},t\in\mathbb{R}.\label{eq: SA-ODE}
\end{align}

As a result, the state trajectory of the dynamical system described
by the ODE (\ref{eq: SA-ODE}) can asymptotically track the state
trajectory of the iteration (\ref{eq: prps-SA}). The convergence
analysis of the stochastic approximation iteration (\ref{eq: prps-SA})
thus can be obtained by analyzing the asymptotic convergence behavior
of the solution to the ODE (\ref{eq: SA-ODE}). This is formally summarized
in the following lemma.

\begin{lemma}
(\textsl{Global Asymptotic Stability of the Limiting ODE {[}35{]})}\label{Lemma: ODE Asympt-Stability}
If the limiting ODE (\ref{eq: SA-ODE}) has a unique equilibrium point
$\mathbf{P}^{*}$ that is globally asymptotically stable, then $\mathbf{P}_{k}$
obtained by stochastic approximation iteration (\ref{eq: prps-SA})
converges to $\mathbf{P}^{*}$ almost surely, i.e., $\mathrm{Pr}\left(\lim_{k\rightarrow\infty}\mathbf{P}_{k}=\mathbf{P}^{*}\right)=1$.
\end{lemma}

\begin{IEEEproof}
Please see Appendix H.
\end{IEEEproof}

In existing literature, the Lyapunov stability theory is utilized
to establish the global asymptotic stability of the limiting ODE (\ref{eq: SA-ODE})
{[}35{]}. Specifically, if there exists a Lyapunov function $V:\mathbb{S}_{+}^{S}\rightarrow\mathbb{R}_{+}$,
such that $\dot{V}\left(\mathbf{P}\left(t\right)\right)=\left\langle \nabla V\left(\mathbf{P}\left(t\right)\right),f\left(\mathbf{P}\left(t\right)\right)\right\rangle \leq0$
with equality if and only if $\mathbf{P}\left(t\right)=\mathbf{P}^{*}$,
then the limiting ODE (\ref{eq: SA-ODE}) is globally asymptotically
stable with a unique equilibrium at $\mathbf{P}\left(t\right)=\mathbf{P}^{*}$.
However, such an approach does not provide guidelines on how to design
and construct the Lyapunov function $V\left(\mathbf{P}\left(t\right)\right)$.
Moreover, the high nonlinearity of $f\left(\mathbf{P}\left(t\right)\right)$
makes it even more difficult to find such a feasible Lyapunov function
$V\left(\mathbf{P}\left(t\right)\right)$. To address this challenge,
we introduce the following virtual fixed-point process $\left\{ \widetilde{\mathbf{P}}_{k},k\geq0\right\} $:
\begin{align}
\widetilde{\mathbf{P}}_{k+1} & =\widetilde{\mathbf{P}}_{k}+\xi f\left(\widetilde{\mathbf{P}}_{k}\right),\widetilde{\mathbf{P}}_{0}=\mathbf{P}_{0},\forall k\geq0,\label{eq: virtual-fixed-point-process}
\end{align}
where $\xi\in\left(0,1\right)$ is a constant. The following lemma
characterizes the relationship between the state trajectory of the
virtual fixed-point process (\ref{eq: virtual-fixed-point-process})
and the state trajectory of the solution to the limiting ODE (\ref{eq: SA-ODE}).

\begin{lemma}
(\textsl{State Trajectories of the Limiting ODE and the Virtual Fixed-point
Process)}\label{Lemma: Difference State trajectory} Let $t_{k}=k\xi,\forall k\geq0$.
Define a continuous, piece-wise linear matrix-valued function $\overline{\mathbf{P}}\left(t\right),t\geq0,$
by $\overline{\mathbf{P}}\left(t_{k}\right)=\widetilde{\mathbf{P}}_{k}$
with linear interpolation on each interval $\left[t_{k},t_{k+1}\right]$
as
\begin{align}
\overline{\mathbf{P}}\left(t\right) & =\widetilde{\mathbf{P}}_{k}+\xi^{-1}\left(t-t_{k}\right)\left(\widetilde{\mathbf{P}}_{k+1}-\widetilde{\mathbf{P}}_{k}\right).\label{eq: Interpolated-virtual-fixed-point-process}
\end{align}
Let $\mathbf{P}^{l}\left(t\right),t\geq l,$ denote the trajectory
of the limiting ODE (\ref{eq: SA-ODE}) with initial condition $\left.\mathbf{P}^{l}\left(t\right)\right|_{t=l}=\overline{\mathbf{P}}\left(l\right),\forall l\in\mathbb{R}_{+}$.
Then, for any $l>0$ and $L>0$, it follows that
\begin{align}
 & \sup_{t\in\left[0,L\right]}\left\Vert \overline{\mathbf{P}}\left(l+t\right)-\mathbf{P}^{l}\left(l+t\right)\right\Vert =\mathcal{O}\left(\xi\right).\label{eq: difference of state trajectories}
\end{align}
\end{lemma}

\begin{IEEEproof}
Please see Appendix H.
\end{IEEEproof}

Lemma \ref{Lemma: Difference State trajectory} states that the gap
between the state trajectory of the virtual fixed-point process (\ref{eq: virtual-fixed-point-process})
and that of the limiting ODE (\ref{eq: SA-ODE}) is $\mathcal{O}\left(\xi\right)$,
which can be made arbitrarily small by letting $\xi\rightarrow0$.
Therefore, the convergence of the state trajectory of the virtual
fixed-point process (\ref{eq: virtual-fixed-point-process}) under
arbitrary $\xi\in\left(0,1\right)$ implies the convergence of the
state trajectory of the limiting ODE (\ref{eq: SA-ODE}), which in
turn leads to the convergence of the stochastic approximation iteration
(\ref{eq: prps-SA}) according to Lemma \ref{Lemma: ODE Asympt-Stability}.

\subsection{Sufficient Condition for Online Learning Convergence}

Since the learned value function $\widetilde{V}_{k}\left(\mathbf{x}_{k}\right)$
(\ref{eq: prps-Vk-iter}) and control action solution $\mathbf{u}_{k}$
(\ref{eq: prps-u-iter}) in the proposed online learning algorithm
are obtained based on the successive update of $\mathbf{P}_{k}$ in
the stochastic approximation iteration (\ref{eq: prps-SA}), the convergence
analysis for the learned value function $\widetilde{V}_{k}\left(\mathbf{x}_{k}\right)$
and control action solution $\mathbf{u}_{k}$ can be obtained by analyzing
the convergence of $\mathbf{P}_{k}$ in (\ref{eq: prps-SA}).

According to Lemma \ref{Lemma: Difference State trajectory}, the
convergence of $\mathbf{P}_{k}$ in (\ref{eq: prps-SA}) is equivalent
to the convergence of the virtual fixed-point process $\widetilde{\mathbf{P}}_{k}$
in (\ref{eq: virtual-fixed-point-process}) under arbitrary $\xi\in\left(0,1\right)$.
Furthermore, based on the structure of the virtual fixed-point process
$\widetilde{\mathbf{P}}_{k}$ in (\ref{eq: virtual-fixed-point-process}),
if $\widetilde{\mathbf{P}}_{k}$ converges to $\widetilde{\mathbf{P}}^{*}$,
then the limiting convergent point $\widetilde{\mathbf{P}}^{*}$ must
be the root of $f(\mathbf{P})=\mathbf{0}$, i.e., $f\left(\widetilde{\mathbf{P}}^{*}\right)=\mathbf{0}$.
As a result, if $\mathbf{P}_{k}$ in the proposed SA iteration (\ref{eq: prps-SA})
converges, it will also converge to the root of $f(\mathbf{P})=\mathbf{0}$.
The full convergence results are formally summarized in the following
Theorem.

\begin{thm}
(\textsl{Almost Sure Convergence of the Proposed Online Learning Algorithm)}\label{Thm: Convergence results}
If one of the three conditions \textsl{(a.1)}, \textsl{(a.2)} and
\textsl{(a.3)}  in Lemma \ref{lemma-impacts-MIMO Fading} is satisfied,
or the condition (\ref{eq:suff-condition-existence}) in Theorem \ref{Thm: suff-condition-existence-uniqueness}
is satisfied, denote the unique root of $f(\mathbf{P})=\mathbf{0}$
as $\mathbf{P}^{*}$, then
\end{thm}

\begin{itemize}
\item \textsl{Convergence of the Virtual Fixed-point Process}: $\widetilde{\mathbf{P}}_{k}$
in the proposed fixed-point iteration (\ref{eq: virtual-fixed-point-process})
converges to $\mathbf{P}^{*}$ almost surely, i.e., $\mathrm{Pr}\left(\lim_{k\rightarrow\infty}\widetilde{\mathbf{P}}_{k}=\mathbf{P}^{*}\right)=1$.
\item \textsl{Conver}\textit{gence of the SA Iteration}: $\mathbf{P}_{k}$
in the proposed SA iteration (\ref{eq: prps-SA}) converges to $\mathbf{P}^{*}$
almost surely, i.e., $\mathrm{Pr}\left(\lim_{k\rightarrow\infty}\mathbf{P}_{k}=\mathbf{P}^{*}\right)=1$.
\item \textsl{Convergence of the Value Function and Control Action}: The
learned value function $\widetilde{V}_{k}\left(\mathbf{x}_{k}\right)$
in (\ref{eq: prps-Vk-iter}) converges to the optimal value function
$\widetilde{V}\left(\mathbf{x}_{k}\right)=\mathbf{x}_{k}^{T}\mathbf{P}^{*}\mathbf{x}_{k}$
in Theorem \ref{Thm: reduced-state-Bellman-Eq} almost surely, i.e.,
\begin{align}
\mathrm{Pr}\left(\lim_{k\rightarrow\infty}\widetilde{V}_{k}\left(\mathbf{x}_{k}\right)=\widetilde{V}\left(\mathbf{x}_{k}\right)\right) & =1.
\end{align}
Moreover, the learned control action $\mathbf{u}_{k}$ in (\ref{eq: prps-u-iter})
converges to the optimal control action $\mathbf{u}_{k}^{*}=\left(\mathrm{\Omega}^{k}\right)^{*}\left(\mathbf{\mathbf{S}}_{k}\right)$
in Theorem \ref{Thm: reduced-state-Bellman-Eq} almost surely, i.e.,
\begin{align}
\mathrm{Pr}\left(\lim_{k\rightarrow\infty}\mathbf{u}_{k}=\mathbf{u}_{k}^{*}\right) & =1,
\end{align}
where
\begin{align}
\mathbf{u}_{k}^{*}= & \left(\mathrm{\Omega}^{k}\right)^{*}\left(\mathbf{\mathbf{S}}_{k}\right)=-(\delta_{k}\mathbf{H}_{k}^{T}\mathbf{B}^{T}\mathbf{P}^{*}\mathbf{B}\mathbf{H}_{k}\nonumber \\
 & +\delta_{k}\mathbf{H}_{k}^{T}\mathbf{M}\mathbf{H}_{k}+\mathbf{R})^{-1}\mathbf{H}_{k}^{T}\mathbf{B}^{T}\mathbf{P}^{*}\mathbf{A}\mathbf{\mathbf{x}}_{k}.\label{eq:optimal control action-1}
\end{align}
\end{itemize}
\begin{IEEEproof}
Please see Appendix H.
\end{IEEEproof}

\section{Numerical Results}

In this section, we compare the performance of the proposed online
optimal control scheme with the following baselines via numerical
simulations.
\begin{itemize}
\item \textbf{Baseline 1} \emph{(Existing Q-learning-based LQR for Static
Channels {[}27{]}-{[}29{]})}: The remote controller adopts the existing
Q-learning-based LQR solution that is designed for static channels
to generate the control actions. Specifically, the Q-function is given
by $Q\left(\mathbf{x}_{k},\mathbf{u}_{k}\right)=[\mathbf{x}_{k}^{T},\mathbf{u}_{k}^{T}]\mathbf{\mathbf{\Psi}}[\mathbf{x}_{k};\mathbf{u}_{k}]^{T}$,
where $\mathbf{\mathbf{\Psi}}\in\mathbb{R}^{\left(S+N_{t}\right)\times\left(S+N_{t}\right)}$
is the kernel matrix. The remote controller uses the Q-learning method
in {[}27{]}-{[}29{]} to obtain $\mathbf{\mathbf{M}}.$ Based on the
learned kernel $\mathbf{\mathbf{\mathbf{\Psi}}}$, the remote controller
generates the control action $\mathbf{u}_{k}^{*}=\arg\min_{\mathbf{u}_{k}}Q\left(\mathbf{x}_{k},\mathbf{u}_{k}\right)$.
\item \textbf{Baseline 2} \emph{(Brute-force Q-learning-based LQR over Wireless
Channels without State Reduction)}: The remote controller brute-force
applies the existing Q-learning-based LQR approach for closed-loop
control over wireless fading channels without state reduction. Both
the CSI $\mathbf{H}_{k}$ and the controller random access state $\delta_{k}$
are state variables in the Q-function. Specifically, the Q-function
is given by $Q\left(\mathbf{x}_{k},\delta_{k},\mathrm{vec}(\mathbf{H}_{k}),\mathbf{u}_{k}\right)=[\mathbf{x}_{k}^{T},\delta_{k},\mathrm{vec}(\mathbf{H}_{k}){}^{T},\mathbf{u}_{k}^{T}]\mathbf{\mathbf{\mathbf{\Psi}}}[\mathbf{x}_{k};\delta_{k},\mathrm{vec}(\mathbf{H}_{k});\mathbf{u}_{k}]^{T}$,
where $\mathbf{\mathbf{\Psi}}\in\mathbb{R}^{\left(1+S+N_{r}\cdot N_{t}+N_{t}\right)\times\left(1+S+N_{r}\cdot N_{t}+N_{t}\right)}$
and $\mathrm{vec}(\mathbf{H}_{k})$ is the column-wised vectorization
of the CSI $\mathbf{H}_{k}$. The remote controller uses the Q-learning
method in {[}27{]}-{[}29{]} to obtain $\mathbf{\mathbf{\mathbf{\Psi}}}.$
Based on the learned kernel $\mathbf{\mathbf{\mathbf{\Psi}}}$, the
remote controller generates the control action $\mathbf{u}_{k}^{*}=\arg\min_{\mathbf{u}_{k}}Q\left(\mathbf{x}_{k},\delta_{k},\mathrm{vec}(\mathbf{H}_{k}),\mathbf{u}_{k}\right)$.
\item \textbf{Baseline 3} \emph{(Genie-aided Optimal LQR Control)}: The
remote controller adopts the genie-aided optimal control solution.
Specifically, the remote controller is assumed to know the $\mathbf{P}^{*}$
that satisfies the NME (\ref{eq: NME}). The remote controller generates
the optimal control action $\mathbf{u}_{k}^{*}$ according to (\ref{eq:optimal control action}).
\end{itemize}

\subsection{Comparison of the Accuracy of Learned Control Actions}

\begin{figure}[tbh]
\begin{centering}
\includegraphics[clip,width=0.86\columnwidth]{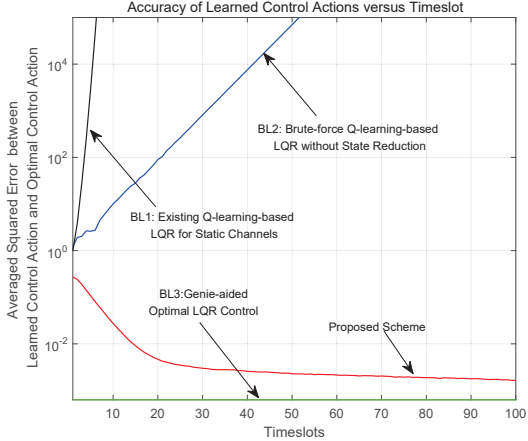}
\par\end{centering}
\caption{\label{fig: uk-error}Comparison of the accuracy of learned control
actions versus the time index $k$. The system parameters are configured
as follows: the state transition matrix $\mathbf{A}=$(0.01, -1.02,
0.3; -0.1, 1.01, 0.2; -0.5, 0.1, 0.2), the control input matrix $\mathbf{B}=(1.1,0.2;-0.2,0.6;-0.3,0.2)$,
the plant noise covariance matrix $\mathbf{W}=0.05\mathbf{I}_{3\times3},$
the LQR weight matrices $\mathbf{Q}=\mathbf{R}=\mathbf{I}_{3\times3}$
and $\mathbf{M}=\mathbf{I}_{2\times2}$, the number of controller
transmit antennas $N_{t}=3$, the number of actuator received antennas
$N_{r}=2$, and the controller random access probability $\mathrm{Pr}\left(\delta_{k}=1\right)=0.5$.}
\end{figure}

Fig. \ref{fig: uk-error} illustrates the accuracy of the learned
control action $\mathbf{u}_{k}$ versus the time index $k$, i.e.,
$\mathbb{E}\left[\left\Vert \mathbf{u}_{k}-\mathbf{u}_{k}^{*}\right\Vert ^{2}\right]$.
It can be observed that with the increase of time, the gap between
the learned control action and the optimal control action becomes
prohibitively large for both Baseline 1 and Baseline 2. This is because
Baseline 1 is designed for static channels. The impacts of random
fading channels and random access of the controller are imprudently
ignored, which causes the divergence of the learned control action
from $\mathbf{u}_{k}^{*}$. Baseline 2 diverges because it is a naive
extension of the standard LQR control solution to wireless channels,
which does not exploit the structure properties of the optimal control
action $\mathbf{u}_{k}^{*}$ w.r.t. the CSI $\mathbf{H}_{k}$ and
random access state $\delta_{k}$. For baseline 3, $\mathbb{E}\left[\left\Vert \mathbf{u}_{k}-\mathbf{u}_{k}^{*}\right\Vert ^{2}\right]=0$
because $\mathbf{P}^{*}$ is assumed known. For the proposed scheme,
the accuracy of the learned control action $\mathbb{E}\left[\left\Vert \mathbf{u}_{k}-\mathbf{u}_{k}^{*}\right\Vert ^{2}\right]$
decreases dramatically as time index $k$ increases. This is because
the control action $\mathbf{u}_{k}$ learned by the proposed scheme
converges almost surely to the optimal control action $\mathbf{u}_{k}^{*}$.

\subsection{Comparison of the Closed-loop Stability}

\begin{figure}[tbh]
\begin{centering}
\includegraphics[clip,width=0.86\columnwidth]{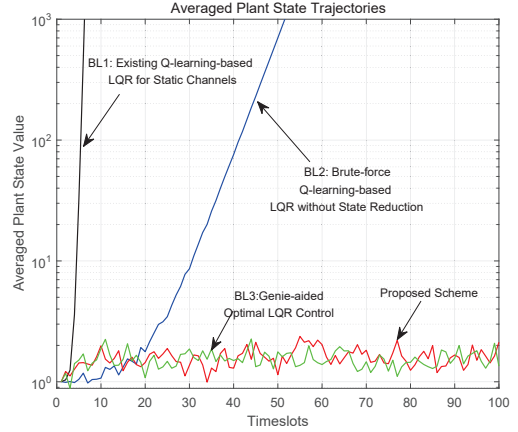}
\par\end{centering}
\caption{\label{fig: average state}Comparison of the closed-loop stability,
i.e., $\mathbb{E}\left[\left\Vert \mathbf{x}_{k}\right\Vert ^{2}\right]$,
versus the time index $k$. The system parameters are configured as
follows: $\mathbf{A}=$(0.01, -1.02, 0.3; -0.1, 1.01, 0.2; -0.5, 0.1,
0.2), $\mathbf{B}=(1.1,0.2;-0.2,0.6;-0.3,0.2)$, $\mathbf{W}=0.05\mathbf{I}_{3\times3},$
$\mathbf{Q}=\mathbf{R}=\mathbf{I}_{3\times3}$, $\mathbf{M}=\mathbf{I}_{2\times2}$,
$N_{t}=3$, $N_{r}=2$, and $\mathrm{Pr}\left(\delta_{k}=1\right)=0.5$.}
\end{figure}

Fig. \ref{fig: average state} illustrates the average state trajectory,
i.e., $\mathbb{E}\left[\left\Vert \mathbf{x}_{k}\right\Vert ^{2}\right]$,
versus the time index $k$. It can be observed that the closed-loop
system is unstable for both Baseline 1 and Baseline 2. This is because
the standard LQR control solution, which is primarily designed for
static channels, fails to achieve closed-loop stability over random
access wireless fading channels. For Baseline 3, the closed-loop system
is stable because the optimal control action is stabilizing. The system
states of the proposed scheme are also stable because of the almost
sure convergence of the learned control action to the stabilizing
optimal control.

\subsection{Comparison of the Computational Complexity}

\begin{figure}[tbh]
\begin{centering}
\includegraphics[clip,width=0.86\columnwidth]{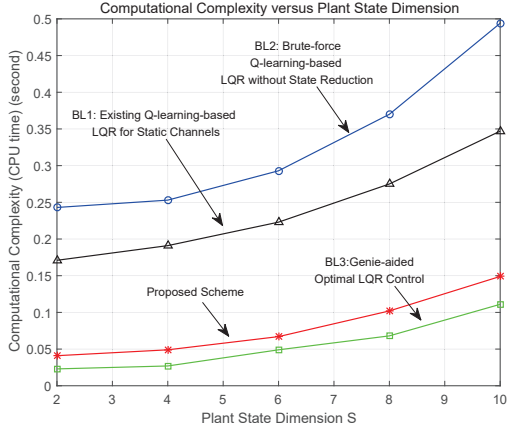}
\par\end{centering}
\caption{\label{fig:cpu time S}Comparison of the CPU time versus the plant
state dimension $S$. The system parameters are configured as follows:
the $\left(i,j\right)$-th element $\mathbf{A}$ is chosen as $(\mathbf{A})_{ij}=-0.1$
when $i=j-1,$ $(\mathbf{A})_{ij}=-0.2$ when $i=j+1,$ $(\mathbf{A})_{ij}=1.01$
when $i=j;$ and $(\mathbf{A})_{ij}=0$ otherwise, with $1\protect\leq i,j\protect\leq S$.
The $\left(i,j\right)$-th element $\mathbf{B}$ is chosen as $(\mathbf{B})_{ij}=(i+j)^{-1}$
with $1\protect\leq i\protect\leq S$ and $1\protect\leq j\protect\leq2$,
$\mathbf{W}=0.05\mathbf{I}_{S\times S},$ $\mathbf{Q}=\mathbf{R}=\mathbf{I}_{S\times S}$,
$\mathbf{M}=\mathbf{I}_{2\times2}$, $N_{t}=3$, $N_{r}=2$, and $\mathrm{Pr}\left(\delta_{k}=1\right)=0.5$.}
\end{figure}

\begin{figure}[tbh]
\begin{centering}
\includegraphics[clip,width=0.86\columnwidth]{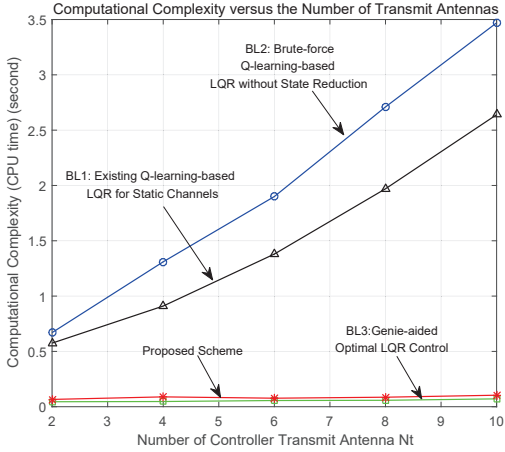}
\par\end{centering}
\caption{\label{fig: CPU time-Nt}Comparison of the CPU time versus the number
of transmit antennas at the remote controller $N_{t}$. The system
parameters are configured as follows: $\mathbf{A}=$(0.01, -1.02,
0.3; -0.1, 1.01, 0.2; -0.5, 0.1, 0.2), $\mathbf{B}=(0.5,1/3;1/3,0.25;0.25,0.2)$,
$\mathbf{W}=0.05\mathbf{I}_{3\times3},$ $\mathbf{Q}=\mathbf{R}=\mathbf{I}_{3\times3}$,
$\mathbf{M}=\mathbf{I}_{2\times2}$, $N_{r}=2$, and $\mathrm{Pr}\left(\delta_{k}=1\right)=0.5$.}
\end{figure}

\begin{figure}[tbh]
\begin{centering}
\includegraphics[clip,width=0.86\columnwidth]{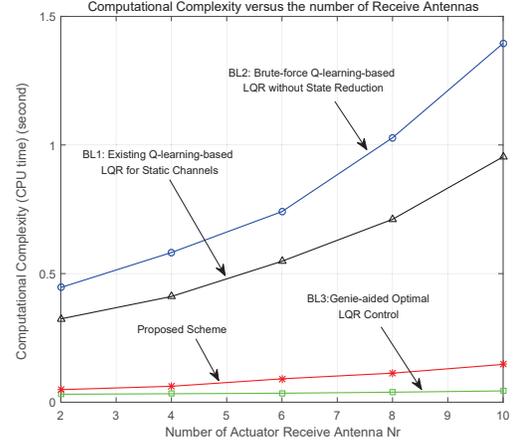}
\par\end{centering}
\caption{\label{fig: CPU time Nr}Comparison of the CPU time versus the number
of receive antennas at the actuator $N_{r}$. The system parameters
are configured as follows: $\mathbf{A}=$(0.01, -1.02, 0.3; -0.1,
1.01, 0.2; -0.5, 0.1, 0.2), the $\left(i,j\right)$-th element $\mathbf{B}$
is chosen as $(\mathbf{B})_{ij}=(i+j)^{-1}$ with $1\protect\leq i\protect\leq3$
and $1\protect\leq j\protect\leq N_{r}$, $\mathbf{W}=0.05\mathbf{I}_{3\times3},$
$\mathbf{Q}=\mathbf{R}=\mathbf{I}_{3\times3}$, $\mathbf{M}=\mathbf{I}_{N_{r}\times N_{r}}$,
$N_{t}=3$, and $\mathrm{Pr}\left(\delta_{k}=1\right)=0.5$.}
\end{figure}

Fig. \ref{fig:cpu time S}, Fig. \ref{fig: CPU time-Nt} and Fig.
\ref{fig: CPU time Nr} illustrates the computational complexity (the
CPU time for $10^{4}$ simulation runs) versus the plant state dimension
$S$, the number of controller transmit antennas $N_{t}$, and the
number of actuator receive antennas $N_{r}$, respectively. It can
be observed that Baseline 2 has the highest computational complexity
because the Q-function to be learned has $\left(1+S+N_{r}\cdot N_{t}+N_{t}\right)^{2}$
dimensions, which is also prohibitively large. Compared with Baseline
1 and Baseline 2, the computational complexity of the proposed scheme
is reduced significantly because of the proposed state reduction technique
in Theorem \ref{Thm: reduced-state-Bellman-Eq}. The computational
complexity gap between the proposed scheme and Baseline 3 is due to
the computation of $\mathbf{P}^{*}$, where the proposed scheme needs
to compute $\mathbf{P}^{*}$ using the proposed SA iteration (\ref{eq: prps-SA}).

\section{Conclusion}

In this paper, we considered the online optimal control over the wireless
MIMO fading channels. We formulated the online optimal control design
as an infinite horizon average cost MDP. We proposed a novel\textsl{
}state reduction technique such that the optimality condition is transformed
into a time-invariant\textsl{ }reduced-state Bellman optimality equation.
We provide the closed-form characterizations on the existence and
uniqueness of the optimal control solution via analyzing the reduced-state
Bellman optimality equation. We further propose a novel SA-based online
learning algorithm that can learn the optimal control action on the
fly based on the plant state observations. We derived a closed-form
sufficient condition that guarantees the almost sure convergence of
the proposed SA-based online learning algorithm to the optimal control
solution. The proposed scheme is also compared with various baselines,
and we show that significant performance gains can be achieved.

\appendix

\subsection{\label{subsec:Proof-of-Thm1-Thm2}Proof of Theorem \ref{Thm: std-Bellman-eq}
and Theorem \ref{Thm: reduced-state-Bellman-Eq}}

Note that if there exists a pair of $\left(\theta,V\left(\mathbf{S}_{k}\right)\right)$
such that the Bellman optimality equation (\ref{eq: standard Bellman Eq})
is satisfied, it follows that for all $\left(\mathbf{\mathbf{S}}_{k},\mathbf{u}_{k}\right)$
pair, the following inequality (\ref{eq:Thm1-inter-1} ) holds:
\begin{align}
 & r\left(\mathbf{\mathbf{S}}_{k},\mathbf{u}_{k}\right)+\sum_{\mathbf{S}_{k+1}}\mathrm{Pr}\left[\left.\mathbf{\mathbf{S}}_{k+1}\right|\mathbf{\mathbf{S}}_{k},\mathbf{u}_{k}\right]V\left(\mathbf{S}_{k+1}\right)\nonumber \\
 & \geq\theta+V\left(\mathbf{S}_{k}\right),\label{eq:Thm1-inter-1}
\end{align}
where the equality sign holds if and only if $\mathbf{u}_{k}=\mathbf{u}_{k}^{*}$
is the minimizer of the L.H.S. of (\ref{eq:Thm1-inter-1}).

Further note that $\sum_{\mathbf{S}_{k+1}}\mathrm{Pr}\left[\left.\mathbf{\mathbf{S}}_{k+1}\right|\mathbf{\mathbf{S}}_{k},\mathbf{u}_{k}\right]V\left(\mathbf{S}_{k+1}\right)=\mathbb{E}\left[\left.V\left(\mathbf{S}_{k+1}\right)\right|\mathbf{\mathbf{S}}_{k},\mathbf{u}_{k}\right]$,
taking full expectation on both sides of (\ref{eq:Thm1-inter-1})
and noting that $\theta$ is a constant, it follows that
\begin{align}
\mathbb{E}\left[r\left(\mathbf{\mathbf{S}}_{k},\mathbf{u}_{k}\right)\right] & +\mathbb{E}\left[V\left(\mathbf{S}_{k+1}\right)\right]\geq\theta+\mathbb{E}\left[V\left(\mathbf{S}_{k}\right)\right].\label{eq:Thm1-inter-2}
\end{align}

Summing the above inequality (\ref{eq:Thm1-inter-2}) of both sides
from $k=0$ to $k=K$, and then dividing both sides by $K$, we have
\begin{align}
\theta & \leq\frac{1}{K}\sum_{k=0}^{K}\mathbb{E}\left[r\left(\mathbf{\mathbf{S}}_{k},\mathbf{u}_{k}\right)\right]+\frac{1}{K}\left(\mathbb{E}\left[V\left(\mathbf{S}_{k+1}\right)-V\left(\mathbf{S}_{0}\right)\right]\right).\label{eq:Thm1-inter-4}
\end{align}
Moreover, based on (\ref{eq:Thm1-inter-1}), for any given $\mathbf{\mathbf{S}}_{k}$,
it follows that
\begin{align}
V\left(\mathbf{S}_{k}\right) & =\limsup_{T\rightarrow\infty}\sum_{t=0}^{T}\mathbb{E}_{\mathbf{u}_{k+t}^{*}}\left[r\left(\mathbf{\mathbf{S}}_{k+t},\mathbf{u}_{k+t}^{*}\right)-\theta\right].\label{eq:Thm1-inter-6}
\end{align}
Substituting (\ref{eq:Thm1-inter-6}) into (\ref{eq:Thm1-inter-4}),
it follows that
\begin{align}
\left(K+T\right)\theta & =\limsup_{\left(K+T\right)\rightarrow\infty}\sum_{k=0}^{K+T}\mathbb{E}\left[r\left(\mathbf{\mathbf{S}}_{k},\mathbf{u}_{k}^{*}\right)\right].\label{eq:Thm1-inter-7}
\end{align}
Therefore,
\begin{align}
\theta & =\limsup_{K\rightarrow\infty}\frac{1}{K}\sum_{k=0}^{K}\mathbb{E}\left[r\left(\mathbf{\mathbf{S}}_{k},\mathbf{u}_{k}^{*}\right)\right].\label{eq:Thm1-inter-8}
\end{align}

As a result, if there exists a pair of $\left(\theta,V\left(\mathbf{S}_{k}\right)\right)$
that satisfies the Bellman optimality equation (\ref{eq: standard Bellman Eq}),
$\theta$, which is independent of any extended states $\mathbf{S}_{0}$,
is the optimal average cost for Problem \ref{problem-1}, and is given
by equation (\ref{eq:Thm1-inter-8}). $V\left(\mathbf{S}_{k}\right)$
is the optimal value function for the extended state $\mathbf{S}_{k}$,
and is given by equation (\ref{eq:Thm1-inter-6}). The optimal control
action $\mathbf{u}_{k}^{*}$ is the minimizer of the R.H.S. of (\ref{eq: reduced-State-Bellman Eq}).
Therefore, Theorem \ref{Thm: std-Bellman-eq} is proved.

Exploiting the i.i.d. property of the MIMO fading channel $\mathbf{H}_{k}$
and the controller random access $\delta_{k}$, the optimality condition
of Problem \ref{problem-1} in Theorem \ref{Thm: std-Bellman-eq}
can be represented as
\begin{align}
 & \theta+V\left(\mathbf{x}_{k},\mathbf{H}_{k},\delta_{k}\right)=\nonumber \\
 & \min_{\mathbf{u}\left(\mathbf{x}_{k},\mathbf{H}_{k},\delta_{k}\right)}\bigg[\widetilde{r}\left(\mathbf{x}_{k},\mathbf{H}_{k},\delta_{k},\mathbf{u}\left(\mathbf{x}_{k},\mathbf{H}_{k},\delta_{k}\right)\right)+\sum_{\mathbf{x}_{k+1},\mathbf{H}_{k+1},\delta_{k+1}}\nonumber \\
 & \mathrm{Pr}\left[\left.\mathbf{x}_{k+1},\mathbf{H}_{k+1},\delta_{k+1}\right|\mathbf{x}_{k},\mathbf{H}_{k},\delta_{k},\mathbf{u}\left(\mathbf{x}_{k},\mathbf{H}_{k},\delta_{k}\right)\right]\nonumber \\
 & \cdot V\left(\mathbf{x}_{k+1},\mathbf{H}_{k+1},\delta_{k+1}\right)\bigg]\nonumber \\
 & =\min_{\mathbf{u}\left(\mathbf{x}_{k},\mathbf{H}_{k},\delta_{k}\right)}\bigg[\widetilde{r}\left(\mathbf{x}_{k},\mathbf{H}_{k},\delta_{k},\mathbf{u}\left(\mathbf{x}_{k},\mathbf{H}_{k},\delta_{k}\right)\right)\nonumber \\
 & +\sum_{\mathbf{x}_{k+1}}\mathrm{Pr}\left[\left.\mathbf{x}_{k+1}\right|\mathbf{x}_{k},\mathbf{H}_{k},\delta_{k},\mathbf{u}\left(\mathbf{x}_{k},\mathbf{H}_{k},\delta_{k}\right)\right]\nonumber \\
 & \cdot(\sum_{\mathbf{H}_{k+1}}\mathrm{Pr}\left[\mathbf{H}_{k+1}\right]\sum_{\delta_{k+1}}\mathrm{Pr}\left[\delta_{k+1}\right]V\left(\mathbf{x}_{k+1},\mathbf{H}_{k+1},\delta_{k+1}\right))\bigg]\nonumber \\
 & =\min_{\mathbf{u}\left(\mathbf{x}_{k},\mathbf{H}_{k},\delta_{k}\right)}\bigg[\widetilde{r}\left(\mathbf{x}_{k},\mathbf{H}_{k},\delta_{k},\mathbf{u}\left(\mathbf{x}_{k},\mathbf{H}_{k},\delta_{k}\right)\right)\nonumber \\
 & \sum_{\mathbf{x}_{k+1}}\mathrm{Pr}\left[\left.\mathbf{x}_{k+1}\right|\mathbf{x}_{k},\mathbf{H}_{k},\delta_{k},\mathbf{u}\left(\mathbf{x}_{k},\mathbf{H}_{k},\delta_{k}\right)\right]\widetilde{V}\left(\mathbf{x}_{k+1}\right)\bigg],\label{eq:Thm2-1}
\end{align}
where $\widetilde{V}\left(\mathbf{x}_{k}\right)=\mathbb{E}\left[\left.V\left(\mathbf{x}_{k+1},\mathbf{H}_{k+1},\delta_{k+1}\right)\right|\mathbf{x}_{k}\right]$.

Taking the conditional expectation (conditioned on $\mathbf{x}_{k}$)
on both sides of (\ref{eq:Thm2-1}), it follows that the reduced state
Bellman optimality equation is given by
\begin{align}
 & \theta+\widetilde{V}\left(\mathbf{x}_{k}\right)=\mathbb{E}\bigg[\min_{\mathbf{u}\left(\mathbf{x}_{k},\mathbf{H}_{k},\delta_{k}\right)}\bigg[\widetilde{r}\left(\mathbf{x}_{k},\mathbf{H}_{k},\delta_{k},\mathbf{u}\left(\mathbf{x}_{k},\mathbf{H}_{k},\delta_{k}\right)\right)\nonumber \\
 & +\sum_{\mathbf{x}_{k+1}}\mathrm{Pr}\left[\left.\mathbf{x}_{k+1}\right|\mathbf{x}_{k},\mathbf{H}_{k},\delta_{k},\mathbf{u}\left(\mathbf{x}_{k},\mathbf{H}_{k},\delta_{k}\right)\right]\widetilde{V}\left(\mathbf{x}_{k+1}\right)\bigg],\forall\mathbf{x}_{k}.\label{eq:Thm2-2}
\end{align}

As a result, if there exists a pair of $\left(\widetilde{\theta},\widetilde{V}\left(\mathbf{x}_{k}\right)\right)$
that solves (\ref{eq:Thm2-2}), then $\widetilde{\theta}=\theta$
is the optimal average cost for Problem \ref{problem-1}, and is given
by equation (\ref{eq:Thm1-inter-8}). $\widetilde{V}\left(\mathbf{x}_{k}\right)=\mathbb{E}\left[\left.V\left(\mathbf{S}_{k}\right)\right|\mathbf{x}_{k}\right]$
is the optimal reduced state value function. The optimal control policy
for Problem \ref{problem-1} is given by $\mathbf{u}^{*}\left(\mathbf{x}_{k},\mathbf{H}_{k},\delta_{k}\right)$,
which attains the minimum of the R.H.S. of (\ref{eq:Thm2-2}). Therefore,
Theorem \ref{Thm: reduced-state-Bellman-Eq} is proved.

\subsection{\label{subsec:Proof-of-Thm-3-Lemma-1}Proof of Theorem \ref{Thm: structure-property-reduced-state-Bellman}
and Lemma \ref{lemma-existence-uniqueness}}

To solve the reduced state Bellman optimality equation (\ref{eq: reduced-State-Bellman Eq}),
we first assume that the reduced state value function $\widetilde{V}\left(\mathbf{x}_{k}\right)$
has a quadratic form of $\mathbf{x}_{k}$ and is given by $\widetilde{V}\left(\mathbf{x}_{k}\right)=\mathbf{x}_{k}^{T}\mathbf{P}\mathbf{x}_{k}$
with $\mathbf{P}\in\mathbb{S}_{+}^{S}$ being a constant positive
definite matrix. Then, equation (\ref{eq: reduced-State-Bellman Eq})
can be represented as
\begin{align}
 & \widetilde{\theta}+\mathbf{x}_{k}^{T}\mathbf{P}\mathbf{x}_{k}=\mathbb{E}\bigg[\min_{\mathbf{u}_{k}}\bigg[\mathbf{x}_{k}^{T}\mathbf{Q}\mathbf{x}_{k}+\mathbf{u}_{k}^{T}\left(\mathbf{R}+\delta_{k}\mathbf{H}_{k}^{T}\mathbf{M}\mathbf{H}_{k}\right)\mathbf{u}_{k}\nonumber \\
 & +\mathrm{Tr}\left(\mathbf{M}\right)+\left(\mathbf{A}\mathbf{x}_{k}+\delta_{k}\mathbf{B}\mathbf{H}_{k}\mathbf{u}_{k}\right)^{T}\mathbf{P}\left(\mathbf{A}\mathbf{x}_{k}+\delta_{k}\mathbf{B}\mathbf{H}_{k}\mathbf{u}_{k}\right)\nonumber \\
 & +\mathrm{Tr}\left(\mathbf{P}\mathbf{W}\right)+\mathrm{Tr}\left(\mathbf{B}^{T}\mathbf{P}\mathbf{B}\right)\bigg],\nonumber \\
 & =\mathbb{E}\left[\min_{\mathbf{u}_{k}}\mathbf{\Phi}\left(\mathbf{x}_{k},\mathbf{u}_{k}\right)\right]+\mathrm{Tr}\left(\mathbf{M}+\mathbf{P}\mathbf{W}+\mathbf{B}^{T}\mathbf{P}\mathbf{B}\right),\forall\mathbf{x}_{k},\label{eq:Thm3-1}
\end{align}
where $\mathbf{\Phi}\left(\mathbf{x}_{k},\mathbf{u}_{k}\right)\in\mathbb{R}_{+}$
is given by
\begin{align}
 & \mathbf{\Phi}\left(\mathbf{x}_{k},\mathbf{u}_{k}\right)=\mathbf{x}_{k}^{T}\mathbf{Q}\mathbf{x}_{k}+\mathbf{u}_{k}^{T}\left(\mathbf{R}+\delta_{k}\mathbf{H}_{k}^{T}\mathbf{M}\mathbf{H}_{k}\right)\mathbf{u}_{k}\nonumber \\
 & +\left(\mathbf{A}\mathbf{x}_{k}+\delta_{k}\mathbf{B}\mathbf{H}_{k}\mathbf{u}_{k}\right)^{T}\mathbf{P}\left(\mathbf{A}\mathbf{x}_{k}+\delta_{k}\mathbf{B}\mathbf{H}_{k}\mathbf{u}_{k}\right)\nonumber \\
 & =\left[\begin{array}{c}
\mathbf{x}_{k}\\
\mathbf{u}_{k}
\end{array}\right]^{T}\mathbf{S}\left[\begin{array}{c}
\mathbf{x}_{k}\\
\mathbf{u}_{k}
\end{array}\right],
\end{align}
and $\mathbf{S}\in\mathbb{R}^{\left(S+N_{r}\right)\times\left(S+N_{r}\right)}$
is given by
\begin{align}
\mathbf{S}= & \left[\begin{array}{cc}
\mathbf{Q}+\mathbf{A}^{T}\mathbf{P}\mathbf{A}, & \delta_{k}\mathbf{A}^{T}\mathbf{P}\mathbf{B}\mathbf{H}_{k}\\
\delta_{k}\mathbf{H}_{k}^{T}\mathbf{B}^{T}\mathbf{P}\mathbf{A}, & \mathbf{R}+\delta_{k}\mathbf{H}_{k}^{T}\mathbf{M}\mathbf{H}_{k}+\delta_{k}\mathbf{H}_{k}^{T}\mathbf{B}^{T}\mathbf{P}\mathbf{B}\mathbf{H}_{k}
\end{array}\right].
\end{align}

Note that the Schur complement of $\mathbf{S}$ is given by
\begin{align}
\mathbf{S}^{c} & =\mathbf{Q}+\mathbf{A}^{T}\mathbf{P}\mathbf{A}-\delta_{k}\mathbf{A}^{T}\mathbf{P}\mathbf{B}\mathbf{H}_{k}(\mathbf{R}+\delta_{k}\mathbf{H}_{k}^{T}\mathbf{M}\mathbf{H}_{k}\nonumber \\
 & +\delta_{k}\mathbf{H}_{k}^{T}\mathbf{B}^{T}\mathbf{P}\mathbf{B}\mathbf{H}_{k})^{-1}\mathbf{H}_{k}^{T}\mathbf{B}^{T}\mathbf{P}\mathbf{A}.
\end{align}
It follows that
\begin{align}
 & \min_{\mathbf{u}_{k}}\mathbf{\Phi}\left(\mathbf{x}_{k},\mathbf{u}_{k}\right)=\mathbf{x}_{k}^{T}\mathbf{S}^{c}\mathbf{x}_{k},\label{eq: Phi-uk}
\end{align}
and the $\mathbf{u}_{k}^{*}$ that achieves the minimum value of (\ref{eq: Phi-uk})
is given by
\begin{align}
 & \mathbf{u}_{k}^{*}\nonumber \\
 & =-\left(\mathbf{R}+\delta_{k}\mathbf{H}_{k}^{T}\mathbf{M}\mathbf{H}_{k}+\delta_{k}\mathbf{H}_{k}^{T}\mathbf{B}^{T}\mathbf{P}\mathbf{B}\mathbf{H}_{k}\right)^{-1}\mathbf{H}_{k}^{T}\mathbf{B}^{T}\mathbf{P}\mathbf{A}\mathbf{x}_{k}.\label{eq:u*-k}
\end{align}

Substituting (\ref{eq: Phi-uk}) and (\ref{eq:u*-k}) into (\ref{eq:Thm3-1}),
the reduced state Bellman optimality equation (\ref{eq: reduced-State-Bellman Eq})
can be represented as
\begin{align}
\widetilde{\theta}+\mathbf{x}_{k}^{T}\mathbf{P}\mathbf{x}_{k} & =\mathrm{Tr}\left(\mathbf{M}+\mathbf{P}\mathbf{W}+\mathbf{B}^{T}\mathbf{P}\mathbf{B}\right)+\mathbf{x}_{k}^{T}\mathbb{E}\left[\mathbf{S}^{c}\right]\mathbf{x}_{k},\forall\mathbf{x}_{k}.\label{eq:Thm3-2}
\end{align}

Assuming $\widetilde{V}\left(\mathbf{x}_{k}\right)$ exists, i.e.,
$\mathbf{P}$ exists, it follows that
\begin{align}
 & \widetilde{\theta}=\mathrm{Tr}\left(\mathbf{M}+\mathbf{P}\mathbf{W}+\mathbf{B}^{T}\mathbf{P}\mathbf{B}\right),\label{eq:Thm-3-theta-tilde}\\
 & \widetilde{V}\left(\mathbf{x}_{k}\right)=\mathbf{x}_{k}^{T}\mathbf{P}\mathbf{x}_{k}=\mathbf{x}_{k}^{T}\mathbb{E}\left[\mathbf{S}^{c}\right]\mathbf{x}_{k}.\forall\mathbf{x}_{k}.\label{eq:Thm3-3}
\end{align}
Therefore, Theorem \ref{Thm: structure-property-reduced-state-Bellman}
is proved.

Note that for any given realization of $\left(\delta_{k},\mathbf{H}_{k}\right)$,
$\widetilde{\theta}$ in (\ref{eq:Thm-3-theta-tilde}), $\widetilde{V}\left(\mathbf{x}_{k}\right)$
in (\ref{eq:Thm3-3}), and $\mathbf{u}_{k}^{*}$ in (\ref{eq:u*-k})
are all uniquely determined by $\mathbf{P}$. Further note that equation
(\ref{eq:Thm3-3}) is satisfied for all $\mathbf{x}_{k}$, it follows
that $\mathbf{P}$ satisfies the following NME:
\begin{align}
\mathbf{P}= & \mathbb{E}\left[\mathbf{S}^{c}\right]=\mathbf{A}^{T}\mathbf{P}\mathbf{A}-\mathbb{E}[\delta_{k}\mathbf{A}^{T}\mathbf{P}\mathbf{B}\mathbf{H}_{k}(\delta_{k}\mathbf{H}_{k}^{T}\mathbf{B}^{T}\mathbf{P}\mathbf{B}\mathbf{H}_{k}\nonumber \\
 & +\delta_{k}\mathbf{H}_{k}^{T}\mathbf{M}\mathbf{H}_{k}+\mathbf{R})^{-1}\mathbf{H}_{k}^{T}\mathbf{B}^{T}\mathbf{P}\mathbf{A}]+\mathbf{Q}.\label{eq: NME-1}
\end{align}

As a result, if the solution $\mathbf{P}^{*}$ to the NME (\ref{eq: NME-1})
is unique, then $\widetilde{\theta}$ , $\widetilde{V}\left(\mathbf{x}_{k}\right)$,
and $\mathbf{u}_{k}^{*}$ all are also unique. Therefore, Lemma \ref{lemma-existence-uniqueness}
is proved.

\subsection{\label{subsec:Proof-of-Lemma-2}Proof of Lemma \ref{lemma-impacts-MIMO Fading}}

According to the Popov--Belevitch--Hautus (PBH) test, the pair $\left(\mathbf{A},\mathbf{B}\right)$
is controllable if and only if there exists no left eigenvector of
$\mathbf{A}$ orthogonal to the columns of $\mathbf{B}$. This means
that provided $\left(\mathbf{A},\mathbf{B}\right)$ is controllable,
if there is a vector-scalar pair $\left(\lambda,\mathbf{v}\right)$,$\lambda\in\mathbb{R}$,
$\mathbf{v}\in\mathbb{R}^{S}$, such that $\mathbf{A}\mathbf{v}=\lambda\mathbf{v}$
and $\mathbf{B}^{T}\mathbf{v}=\mathbf{0},$ then $\mathbf{v}=\mathbf{0}.$

We first prove the following proposition.
\begin{prop}
\label{Prps-1} $\delta_{k}\mathbf{B}\mathbf{H}_{k}$ is statistically
identical to \textbf{$\delta_{k}\left(\mathbf{B}\mathbf{B}^{T}\right)^{\frac{1}{2}}\widetilde{\mathbf{H}}_{k}$},
where each element of $\widetilde{\mathbf{H}}_{k}\in\mathbb{R}^{S\times N_{t}}$
is i.i.d. Gaussian distributed with zero mean and unit variance.
\end{prop}

\begin{IEEEproof}
For any realization of $\delta_{k}$, if $\mathrm{Rank}\left(\mathbf{B}\mathbf{B}^{T}\right)=S$,
i.e., $\left(\mathbf{B}\mathbf{B}^{T}\right)^{\frac{1}{2}}$ is full
rank, then we choose $\widetilde{\mathbf{H}}_{k}$ to be
\begin{align}
\widetilde{\mathbf{H}}_{k}= & \left(\mathbf{B}\mathbf{B}^{T}\right)^{-\frac{1}{2}}\mathbf{B}\mathbf{H}_{k}.
\end{align}
Since each element of $\mathbf{H}_{k}$ is i.i.d. Gaussian distributed
with zero mean and unit variance, it follows that each element of
$\widetilde{\mathbf{H}}_{k}$ is also i.i.d. Gaussian distributed
with zero mean and unit variance.

In the case that $\left(\mathbf{B}\mathbf{B}^{T}\right)^{\frac{1}{2}}$
is rank deficient, let $\mathrm{Rank}\left(\left(\mathbf{B}\mathbf{B}^{T}\right)^{\frac{1}{2}}\right)=\eta_{B}<S$.
Denote the singular value decomposition of $\left(\mathbf{B}\mathbf{B}^{T}\right)^{\frac{1}{2}}$
be $\left(\mathbf{B}\mathbf{B}^{T}\right)^{\frac{1}{2}}=\mathbf{U}^{T}\left(\boldsymbol{\Lambda}\right)^{\frac{1}{2}}\mathbf{U}$,
where $\boldsymbol{\Lambda}=\mathrm{diag}\left(\lambda_{1},\cdots,\lambda_{\eta_{B}},0,\cdots,0\right)$,
$\lambda_{i}$ ($1\leq i\leq\eta_{B}$) are the $\eta_{B}$ nonzero
singular values of $\left(\mathbf{B}\mathbf{B}^{T}\right)$. Denote
$\widetilde{\mathbf{\boldsymbol{\Lambda}}}=\mathrm{diag}\left(\left(\lambda_{1}\right)^{-1},\cdots,\left(\lambda_{\eta_{B}}\right)^{-1},0,\cdots,0\right)$,
denote $\left[\left(\widetilde{\mathbf{H}}\left(k\right)\right)_{l,m}\right]$
as the $S\times S$ dimensional matrix with all the elements being
0 except the $l$-th row and $m$-th column element being $\left(\widetilde{\mathbf{H}}\left(k\right)\right)_{l,m}$.
Denote $h_{l,m}\left(k\right)=\left(\mathbf{U}_{2}^{T}\widetilde{\mathbf{\boldsymbol{\Lambda}}}\mathbf{U}_{2}\mathbf{B}\mathbf{H}\left(k\right)\right)_{l,m}$.
Let $g_{l,m}\left(k\right)$ an i.i.d. Gaussian distributed random
variable with zero mean and unit variance. In this case, $\widetilde{\mathbf{H}}\left(k\right)$
is given by:
\begin{align}
\left(\widetilde{\mathbf{H}}\left(k\right)\right)_{l,m} & =\begin{cases}
h_{l,m}\left(k\right), & if\ h_{l,m}\left(k\right)\neq0;\\
g_{l,m}\left(k\right), & otherwise.
\end{cases}
\end{align}
Therefore, Proposition \ref{Prps-1} is proved.
\end{IEEEproof}

In the following, we prove Lemma \ref{lemma-impacts-MIMO Fading}
based on the above PBH test and Proposition 1.

\textsl{Proof of (a.1)}: Given $\delta_{k}=1$, suppose that a there
is a vector-scalar pair $\left(\widetilde{\lambda},\widetilde{\mathbf{v}}\right)$
such that $\mathbf{A}\widetilde{\mathbf{v}}=\widetilde{\lambda}\widetilde{\mathbf{v}}$
and 
\begin{align}
 & \left(\left(\mathbf{B}\mathbf{B}^{T}\right)^{\frac{1}{2}}\widetilde{\mathbf{H}}_{k}\right)^{T}\widetilde{\mathbf{v}}=\left(\widetilde{\mathbf{H}}_{k}\right)^{T}\left(\mathbf{B}\mathbf{B}^{T}\right)^{\frac{1}{2}}\widetilde{\mathbf{v}}=\mathbf{0}.\label{eq: lemma2-1}
\end{align}
Since $N_{t}\geq S$, it follows that $\widetilde{\mathbf{H}}_{k}\left(\widetilde{\mathbf{H}}_{k}\right)^{T}$
is full rank w.p.1.. Multiplying $\left(\widetilde{\mathbf{H}}_{k}\left(\widetilde{\mathbf{H}}_{k}\right)^{T}\right)^{-1}\widetilde{\mathbf{H}}_{k}$
on both sides of (\ref{eq: lemma2-1}), it follows that $\left(\mathbf{B}\mathbf{B}^{T}\right)^{\frac{1}{2}}\widetilde{\mathbf{v}}=\mathbf{0}.$
Since $\left(\mathbf{B}\mathbf{B}^{T}\right)^{\frac{1}{2}}$ and $\left(\mathbf{B}\mathbf{B}^{T}\right)$
have the same null-space, we conclude that 
\begin{align}
 & \mathbf{B}\mathbf{B}^{T}\widetilde{\mathbf{v}}=\mathbf{0}.\label{eq:lemma2-2}
\end{align}

Multiplying $\widetilde{\mathbf{v}}^{T}$ on both sides of (\ref{eq:lemma2-2}),
it follows that $\left\Vert \mathbf{B}^{T}\widetilde{\mathbf{v}}\right\Vert _{2}^{2}=0$,
which leads to $\mathbf{B}^{T}\widetilde{\mathbf{v}}=\mathbf{0}$.
Since the pair $\left(\mathbf{A},\mathbf{B}\right)$ is controllable,
we conclude that $\widetilde{\mathbf{v}}=\mathbf{0}$. As a result,
$\left(\mathbf{A},\mathbf{B}\mathbf{H}_{k}\right)$ is almost surely
controllable. Therefore, \textsl{(a.1)} in Lemma \ref{lemma-impacts-MIMO Fading}
is proved.

\textsl{Proof of (b.1)}: Based on the proof of \textsl{(a.1)}, we
know that given $N_{t}\geq S$ and $\delta_{k}=1$, $\left(\mathbf{A},\mathbf{B}\mathbf{H}_{k}\right)$
is almost surely controllable. Therefore, in the case $N_{r}\geq S$
and $0<\mathrm{Pr}\left(\delta_{k}=1\right)<1$, $\left(\mathbf{A},\delta_{k}\mathbf{B}\mathbf{H}_{k}\right)$
is almost surely controllable when $\delta_{k}=1$, and when $\delta_{k}=0$,
$\left(\mathbf{A},\mathbf{0}\right)$ is uncontrollable. Therefore,
\textsl{(b.1)} in Lemma \ref{lemma-impacts-MIMO Fading} is proved.

In the case that $N_{t}<S$, $\widetilde{\mathbf{H}}_{k}\left(\widetilde{\mathbf{H}}_{k}\right)^{T}$
is rank deficient and $\mathrm{Rank}\left(\widetilde{\mathbf{H}}_{k}\right)=N_{r}$
w.p.1.. We have 
\begin{align}
 & \left(\widetilde{\mathbf{H}}_{k}\right)^{T}\mathbf{B}\mathbf{B}^{T}\widetilde{\mathbf{v}}=\left(\widetilde{\mathbf{H}}_{k}\right)^{T}\mathbf{U}^{T}\boldsymbol{\Xi}\mathbf{U}\widetilde{\mathbf{v}}=\widehat{\mathbf{H}}_{k}\boldsymbol{\Lambda}\widehat{\mathbf{v}}=\nonumber \\
 & \left[\left(\widehat{\mathbf{H}}_{k}\right)_{1:N_{t};1:\eta_{B}}\cdot\mathrm{diag}\left(\lambda_{1},\cdots,\lambda_{\eta_{B}}\right),\mathbf{0}_{N_{t}\times\left(S-\eta_{B}\right)}\right]\left[\begin{array}{c}
\widehat{\mathbf{v}}_{1}\\
\widehat{\mathbf{v}}_{2}
\end{array}\right],
\end{align}
where $\widehat{\mathbf{H}}_{k}=\left(\mathbf{U}\widetilde{\mathbf{H}}_{k}\right)^{T}\in\mathbb{R}^{N_{t}\times S}$
is a random matrix with each element being i.i.d. Gaussian distributed
with zero mean and unit variance, $\widehat{\mathbf{v}}_{1}\in\mathbb{R}^{\eta_{B}\times1}$,
$\widehat{\mathbf{v}}_{2}\in\mathbb{R}^{\left(S-\eta_{B}\right)\times1}$
and $\widehat{\mathbf{v}}=\left[\widehat{\mathbf{v}}_{1}^{T},\widehat{\mathbf{v}}_{2}^{T}\right]^{T}=\mathbf{U}\widetilde{\mathbf{v}}$.

Now suppose $\left(\widetilde{\mathbf{H}}_{k}\right)^{T}\mathbf{B}\mathbf{B}^{T}\widetilde{\mathbf{v}}=\mathbf{0}$,
it follows that $\widehat{\mathbf{v}}_{1}$ must lie in the null space
of $\left(\widehat{\mathbf{H}}_{k}\right)_{1:N_{t};1:\eta_{B}}\cdot\mathrm{diag}\left(\lambda_{1},\cdots,\lambda_{\eta_{B}}\right)$,
whereas $\widehat{\mathbf{v}}_{2}$ can take any value provided that
$\widehat{\mathbf{v}}_{2}\in\mathbb{R}^{\left(S-\eta_{B}\right)\times1}$.
In the following, we separate $N_{t}<S$ into two sub-cases: $\eta_{B}\leq N_{t}<S$
and $N_{t}<\eta_{B}$.

\textsl{Proof of (a.2)}: In the case that $\eta_{B}\leq N_{t}<S$,
the left inverse of $\left(\widehat{\mathbf{H}}_{k}\right)_{1:N_{t};1:\eta_{B}}$
exists, and
\begin{align}
 & \left(\widehat{\mathbf{H}}_{k}\right)_{1:N_{t};1:\eta_{B}}\cdot\mathrm{diag}\left(\lambda_{1},\cdots,\lambda_{\eta_{B}}\right)\widehat{\mathbf{v}}_{1}=\mathbf{0}\Longleftrightarrow\widehat{\mathbf{v}}_{1}=\mathbf{0}.\label{eq:lemma-2-3}
\end{align}
Moreover, note that $\mathbf{A}\widetilde{\mathbf{v}}=\mathbf{A}\mathbf{U}^{T}\mathbf{U}\widetilde{\mathbf{v}}$,
it follows that
\begin{align}
 & \mathbf{U}\mathbf{A}\widetilde{\mathbf{v}}=\left(\mathbf{U}\mathbf{A}\mathbf{U}^{T}\right)\mathbf{U}\widetilde{\mathbf{v}}=\widetilde{\mathbf{A}}\widehat{\mathbf{v}}\nonumber \\
 & =\left[\begin{array}{cc}
\widetilde{\mathbf{A}}_{11} & \widetilde{\mathbf{A}}_{12}\\
\widetilde{\mathbf{A}}_{21} & \widetilde{\mathbf{A}}_{22}
\end{array}\right]\left[\begin{array}{c}
\boldsymbol{0}\\
\widehat{\mathbf{v}}_{2}
\end{array}\right]=\left[\begin{array}{c}
\widetilde{\mathbf{A}}_{12}\widehat{\mathbf{v}}_{2}\\
\widetilde{\mathbf{A}}_{22}\widehat{\mathbf{v}}_{2}
\end{array}\right].\label{eq:lemma-2-4}
\end{align}
Therefore, $\widetilde{\mathbf{A}}\widehat{\mathbf{v}}=\lambda\widehat{\mathbf{v}}$
is equivalent to $\widetilde{\mathbf{A}}_{12}\widehat{\mathbf{v}}_{2}=\mathbf{0}$
and $\widetilde{\mathbf{A}}_{22}\widehat{\mathbf{v}}_{2}=\lambda\widehat{\mathbf{v}}_{2}$.
It follows that $\widehat{\mathbf{v}}_{2}=\mathbf{0}$ if and only
if the pair $\left(\widetilde{\mathbf{A}}_{22},\widetilde{\mathbf{A}}_{12}^{T}\right)$
is controllable. As a result, we conclude that when $N_{t}\geq\eta_{B}$
and $\left(\widetilde{\mathbf{A}}_{22},\widetilde{\mathbf{A}}_{12}^{T}\right)$
is controllable, the $\widetilde{\mathbf{v}}$ that simultaneously
satisfies $\left(\widetilde{\mathbf{H}}_{k}\right)^{T}\mathbf{B}\mathbf{B}^{T}\widetilde{\mathbf{v}}=\mathbf{0}$
and $\widetilde{\mathbf{A}}\widehat{\mathbf{v}}=\lambda\widehat{\mathbf{v}}$
is $\widehat{\mathbf{v}}=\mathbf{0}$, i.e., $\left(\mathbf{A},\mathbf{B}\mathbf{H}_{k}\right)$
is controllable w.p.1.. Therefore, \textsl{(a.2)} in Lemma \ref{lemma-impacts-MIMO Fading}
is proved.

\textsl{Proof of (a.3)}: In the case that $N_{t}<\eta_{B}$, the left
inverse of $\left(\widehat{\mathbf{H}}_{k}\right)_{1:N_{t};1:\eta_{B}}$
does not exist and equation (\ref{eq:lemma-2-3}) no longer holds,
and the analysis approach for $N_{t}\geq\eta_{B}$ cannot be applied.
In this case, we shall exploit the null space property of the Gaussian
random matrix {[}36{]}. Based on Lemma 1 in {[}36{]}, it follows that
\begin{align}
 & \left(\widehat{\mathbf{H}}_{k}\right)_{1:N_{t};1:\eta_{B}}\cdot\mathrm{diag}\left(\lambda_{1},\cdots,\lambda_{\eta_{B}}\right)\widehat{\mathbf{v}}_{1}=\mathbf{0}\nonumber \\
 & \Longleftrightarrow\widehat{\mathbf{v}}_{1}=\mathrm{diag}\left(\lambda_{1}^{-1},\cdots,\lambda_{\eta_{B}}^{-1}\right)\mathbf{Z}\mathbf{n},\label{eq:lemma-2-5}
\end{align}
where $\mathbf{Z}\in\mathbb{R}^{\eta_{B}\times\left(\eta_{B}-N_{t}\right)}$
has i.i.d. $\mathcal{N}\left(0,1\right)$ entries and $\mathbf{n}\in\mathbb{R}^{\left(\eta_{B}-N_{t}\right)\times1}$
is a constant vector. $\widehat{\mathbf{v}}_{2}$ still can take any
value provided that $\widehat{\mathbf{v}}_{2}\in\mathbb{R}^{\left(S-\eta_{B}\right)\times1}$.
Now suppose $\left(\mathbf{A},\mathbf{B}\mathbf{H}_{k}\right)$ is
uncontrollable with positive probability, it follows that there are
a total number $\left(\eta_{B}-N_{t}\right)$ pairs of $\left(\widehat{\mathbf{v}}_{1}^{\left(i\right)},\widehat{\mathbf{v}}_{2}^{\left(i\right)}\right)$
such that $\widetilde{\mathbf{A}}\widehat{\mathbf{V}}=\lambda\widehat{\mathbf{V}}$,
where $\widehat{\mathbf{V}}\in\mathbb{R}^{S\times N_{t}}$ is given
by
\begin{align}
 & \widehat{\mathbf{V}}=\left[\begin{array}{ccc}
\widehat{\mathbf{v}}_{1}^{\left(1\right)}, & \cdots, & \widehat{\mathbf{v}}_{1}^{\left(\eta_{B}-N_{t}\right)}\\
\widehat{\mathbf{v}}_{2}^{\left(1\right)}, & \cdots, & \widehat{\mathbf{v}}_{2}^{\left(\eta_{B}-N_{t}\right)}
\end{array}\right],
\end{align}
where $\widehat{\mathbf{v}}_{1}^{\left(1\right)},\cdots,\widehat{\mathbf{v}}_{1}^{\left(\eta_{B}-N_{t}\right)}$
are linearly independent, i.e., $\mathrm{Rank}\left(\widehat{\mathbf{V}}\right)\geq\left(\eta_{B}-N_{t}\right)$.
Therefore, there exists a $\widehat{\mathbf{V}}\neq\mathbf{0}$ such
that $\widetilde{\mathbf{A}}\widehat{\mathbf{V}}=\lambda\widehat{\mathbf{V}}$
if and only if $\mathrm{dim}\left(\mathrm{Null\left(\widetilde{\mathbf{A}}-\lambda\mathbf{I}\right)}\right)\geq\left(\eta_{B}-N_{t}\right)$.
Note that if $\mathrm{Rank}\left(\widetilde{\mathbf{A}}-\lambda\mathbf{I}\right)>\left(S-\eta_{B}+N_{t}\right)$,
$\forall\lambda\in\mathbb{C}$, then $\mathrm{dim}\left(\mathrm{Null\left(\widetilde{\mathbf{A}}-\lambda\mathbf{I}\right)}\right)<\left(\eta_{B}-N_{t}\right)$.
Therefore, when $\mathrm{Rank}\left(\widetilde{\mathbf{A}}-\lambda\mathbf{I}\right)>\left(S-\eta_{B}+N_{t}\right)$,
$\forall\lambda\in\mathbb{C}$, $\widetilde{\mathbf{A}}\widehat{\mathbf{V}}=\lambda\widehat{\mathbf{V}}$
if and only if $\widehat{\mathbf{V}}=\mathbf{0}$. Therefore, we conclude
that when $N_{t}<\eta_{B}$ and $\mathrm{Rank}\left(\widetilde{\mathbf{A}}-\lambda\mathbf{I}\right)>\left(S-\eta_{B}+N_{t}\right),\forall\lambda\in\mathbb{C}$,
the $\widetilde{\mathbf{v}}$ that simultaneously satisfies $\left(\widetilde{\mathbf{H}}_{k}\right)^{T}\mathbf{B}\mathbf{B}^{T}\widetilde{\mathbf{v}}=\mathbf{0}$
and $\widetilde{\mathbf{A}}\widehat{\mathbf{v}}=\lambda\widehat{\mathbf{v}}$
is $\widehat{\mathbf{v}}=\mathbf{0}$, i.e., $\left(\mathbf{A},\mathbf{B}\mathbf{H}_{k}\right)$
is controllable w.p.1.. Therefore, \textsl{(a.3)} in Lemma \ref{lemma-impacts-MIMO Fading}
is proved.

\textsl{Proof of (b.2), (b.3)}: Note that given $\delta_{k}=1$, condition
\textsl{(b.2) and (b.3)} are reduced to condition \textsl{(a.2) and
(a.3)}, respectively\textsl{. }Further note that when $\delta_{k}=0$,
$\left(\mathbf{A},\mathbf{0}\right)$ is uncontrollable, it follows
that \textsl{(b.2), and (b.3)} in Lemma \ref{lemma-impacts-MIMO Fading}
is proved.

\textsl{Proof of (c)}: Based on the proof of case \textsl{(a) and
(b)}, it follows that in the case $\eta_{B}\leq N_{t}<S$, if $\left(\widetilde{\mathbf{A}}_{22},\widetilde{\mathbf{A}}_{12}^{T}\right)$
is uncontrollable, then $\left(\mathbf{A},\delta_{k}\mathbf{B}\mathbf{H}_{k}\right)$
is uncontrollable w.p.1. regardless of the realization of $\delta_{k}$.
In the case that $N_{t}<\eta_{B}$, if $\exists\lambda\in\mathbb{C}$
such that $\mathrm{Rank}\left(\widetilde{\mathbf{A}}-\lambda\mathbf{I}\right)\leq\left(S-\eta_{B}+N_{t}\right),$
then $\left(\mathbf{A},\mathbf{B}\mathbf{H}_{k}\right)$ is uncontrollable
w.p.1. regardless of the realization of $\delta_{k}$. Therefore,
\textsl{(c)} in Lemma \ref{lemma-impacts-MIMO Fading} is proved.

\subsection{\label{subsec:Proof-of-Thm-4-Lemma-3}Proof of Theorem \ref{Thm-PSD decomposition}
and Lemma \ref{lemma-decomposed-NME}}

Applying Lemma 4 in {[}37{]}, it follows that, for given realizations
of $\mathbf{H}_{k}$ and $\delta_{k}$, the NME (\ref{eq: NME}) can
be represented as
\begin{align}
\mathbf{P}= & \mathbf{A}^{T}\mathbf{P}\mathbf{A}-\mathbf{A}^{T}\mathbf{P}\mathbf{\Psi}_{k}(\mathbf{\Psi}_{k}^{T}\mathbf{P}\mathbf{\Psi}_{k}+\mathbf{I})^{-1}\mathbf{\Psi}_{k}^{T}\mathbf{P}\mathbf{A}]+\mathbf{Q}\nonumber \\
= & \mathbf{A}\Big(\mathbf{\Psi}_{k}\mathbf{\Psi}_{k}^{T}+\left(\mathbf{P}\right)^{-1}\Big)^{-1}\mathbf{A}^{T}+\mathbf{Q},\label{eq:Thm-4-1}
\end{align}
where
\begin{align}
 & \mathbf{\Psi}_{k}=\delta_{k}\mathbf{B}\mathbf{H}_{k}\left(\delta_{k}\mathbf{H}_{k}^{T}\mathbf{M}\mathbf{H}_{k}+\mathbf{R}\right)^{-\frac{1}{2}}.
\end{align}
Let the SVD of $\mathbf{\Psi}_{k}\mathbf{\Psi}_{k}^{T}$ be $\mathbf{\Psi}_{k}\mathbf{\Psi}_{k}^{T}=\mathbf{V}_{k}^{T}\boldsymbol{\Pi}_{k}\mathbf{V}_{k}$,
it follows equation (\ref{eq:Thm-4-1}) can be further simplified
as
\begin{align}
 & \mathbf{P}=\mathbf{A}^{T}\mathbf{V}_{k}^{H}\Big(\boldsymbol{\Pi}_{k}+\left(\mathbf{V}_{k}\mathbf{P}\mathbf{V}_{k}^{H}\right)^{-1}\Big)^{-1}\mathbf{V}_{k}\mathbf{A}+\mathbf{Q}.
\end{align}
We now decompose $\mathbf{V}_{k}\mathbf{P}\mathbf{V}_{k}^{H}$ into
two parts with $\mathbf{V}_{k}\mathbf{P}\mathbf{V}_{k}^{H}=\widetilde{\boldsymbol{\mathbf{P}}}_{k}^{uc}+\widetilde{\boldsymbol{\mathbf{P}}}_{k}^{c}$.
Let 
\begin{align}
\widetilde{\boldsymbol{\mathbf{P}}}_{k}^{uc} & =\left(\mathbf{I}-\boldsymbol{\Lambda}_{k}\right)\mathbf{V}_{k}\mathbf{P}\mathbf{V}_{k}^{H}\left(\mathbf{I}-\boldsymbol{\Lambda}_{k}\right)\nonumber \\
 & -\left[\begin{array}{cc}
\mathbf{0} & \mathbf{0}\\
\mathbf{0} & \mathbf{L}_{k}^{T}\left(\mathbf{V}_{k}\mathbf{P}\mathbf{V}_{k}^{H}\right)_{\gamma_{k}}\mathbf{L}_{k}
\end{array}\right],\\
\widetilde{\boldsymbol{\mathbf{P}}}_{k}^{c} & =\left[\begin{array}{cc}
\left(\mathbf{V}_{k}\mathbf{P}\mathbf{V}_{k}^{H}\right)_{\gamma_{k}} & \left(\mathbf{V}_{k}\mathbf{P}\mathbf{V}_{k}^{H}\right)_{\gamma_{k}}\mathbf{L}_{k}\\
\mathbf{L}_{k}^{T}\left(\mathbf{V}_{k}\mathbf{P}\mathbf{V}_{k}^{H}\right)_{\gamma_{k}} & \mathbf{L}_{k}^{T}\left(\mathbf{V}_{k}\mathbf{P}\mathbf{V}_{k}^{H}\right)_{\gamma_{k}}\mathbf{L}_{k}
\end{array}\right].
\end{align}
It is easy to verify that $\widetilde{\boldsymbol{\mathbf{P}}}_{k}^{uc}\in\mathbb{S}_{+}^{S}$
and $\boldsymbol{\Pi}_{k}\widetilde{\boldsymbol{\mathbf{P}}}_{k}^{uc}=\mathbf{0}$;
$\widetilde{\boldsymbol{\mathbf{P}}}_{k}^{c}\in\mathbb{S}_{+}^{S}$
and $\mathrm{\text{ker}}\left(\boldsymbol{\Pi}_{k}\widetilde{\boldsymbol{\mathbf{P}}}_{k}^{c}\right)=\mathrm{\text{ker}}\left(\widetilde{\boldsymbol{\mathbf{P}}}_{k}^{c}\right)$.
Therefore, $\boldsymbol{\mathbf{P}}_{k}^{uc}=\mathbf{V}_{k}^{T}\widetilde{\boldsymbol{\mathbf{P}}}_{k}^{uc}\mathbf{V}_{k}$
and $\boldsymbol{\mathbf{P}}_{k}^{uc}=\mathbf{V}_{k}^{T}\widetilde{\boldsymbol{\mathbf{P}}}_{k}^{uc}\mathbf{V}_{k}$.
Therefore, Theorem \ref{Thm-PSD decomposition} is proved.

Substitute (\ref{eq:Pk-c}) and (\ref{eq:Pk-uc}) into (\ref{eq: NME}),
it follows that
\begin{align}
 & \mathbf{P}=\mathbb{E}[\mathbf{A}^{T}\left(\boldsymbol{\mathbf{P}}_{k}^{uc}+\boldsymbol{\mathbf{P}}_{k}^{c}\right)\mathbf{A}-\mathbf{A}^{T}\left(\boldsymbol{\mathbf{P}}_{k}^{uc}+\boldsymbol{\mathbf{P}}_{k}^{c}\right)\mathbf{\Psi}_{k}\nonumber \\
 & \cdot\Big(\mathbf{\Psi}_{k}\left(\boldsymbol{\mathbf{P}}_{k}^{uc}+\boldsymbol{\mathbf{P}}_{k}^{c}\right)\mathbf{\Psi}_{k}^{T}+\mathbf{I}\Big)^{-1}\mathbf{\Psi}_{k}^{T}\left(\boldsymbol{\mathbf{P}}_{k}^{uc}+\boldsymbol{\mathbf{P}}_{k}^{c}\right)\mathbf{A}]+\mathbf{Q}.\label{eq:Thm-4-2}
\end{align}
Further note that $\mathbf{\Psi}_{k}^{T}\boldsymbol{\mathbf{P}}_{k}^{uc}=\boldsymbol{0}$,
it follows that (\ref{eq:Thm-4-2}) can be simplified as
\begin{align}
 & \mathbf{P}=\mathbb{E}[\mathbf{A}^{T}\left(\boldsymbol{\mathbf{P}}_{k}^{uc}+\boldsymbol{\mathbf{P}}_{k}^{c}\right)\mathbf{A}-\mathbf{A}^{T}\boldsymbol{\mathbf{P}}_{k}^{c}\mathbf{\Psi}_{k}\nonumber \\
 & \cdot\Big(\mathbf{\Psi}_{k}\boldsymbol{\mathbf{P}}_{k}^{c}\mathbf{\Psi}_{k}^{T}+\mathbf{I}\Big)^{-1}\mathbf{\Psi}_{k}^{T}\boldsymbol{\mathbf{P}}_{k}^{c}\mathbf{A}]+\mathbf{Q}.\label{eq:Thm-4-3}
\end{align}
Therefore, Lemma \ref{lemma-decomposed-NME} is proved.

\subsection{\label{subsec:Proof-of-Thm-5}Proof of Theorem \ref{Thm: suff-condition-existence-uniqueness}}

We first prove the existence of $\mathbf{P}^{*}$ that satisfies the
NME (\ref{eq: NME}) when the closed-loop control system is almost
surely controllable, i.e., one of the three conditions \textsl{(a.1)},
\textsl{(a.2)} and \textsl{(a.3)} in Lemma \ref{lemma-impacts-MIMO Fading}
is satisfied. Specifically, denote $g(\mathbf{P})=f(\mathbf{P})+\mathbf{P}$,
it suffices to prove there is a $\mathbf{P}^{*}$ such that $\mathbf{P}^{*}=g\left(\mathbf{P}^{*}\right)$.
Note that there is a $\mathbf{P}^{\left(1\right)}=\mathbf{0}$ such
that $\mathbf{P}^{\left(1\right)}=\mathbf{0}<g\left(\mathbf{P}^{\left(1\right)}\right)=\mathbf{Q}$.
Furthermore, for any given realization of $\mathbf{H}_{k}=\mathbf{H}$,
$g(\left.\mathbf{P}\right|\mathbf{H}_{k}=\mathbf{H})$ can be represented
as

\begin{align}
 & g(\mathbf{P}|\mathbf{H}_{k}=\mathbf{H})=\mathbf{A}^{T}\left(\mathbf{I}-\mathbf{K}\mathbf{\Psi}\right)\mathbf{P}\left(\mathbf{I}-\mathbf{K}\mathbf{\Psi}\right)^{T}\mathbf{A}\nonumber \\
 & +\mathbf{A}^{T}\mathbf{K}\mathbf{K}^{T}\mathbf{A}+\mathbf{Q}\nonumber \\
 & \leq\mathbf{A}^{T}\left(\mathbf{I}-\widetilde{\mathbf{K}}\mathbf{\Psi}\right)\mathbf{P}\left(\mathbf{I}-\widetilde{\mathbf{K}}\mathbf{\Psi}\right)^{T}\mathbf{A}+\mathbf{A}^{T}\widetilde{\mathbf{K}}\mathbf{\widetilde{\mathbf{K}}}^{T}\mathbf{A}+\mathbf{Q},
\end{align}
where $\mathbf{\Psi}=\mathbf{\Psi}_{k}|_{\mathbf{H}_{k}=\mathbf{H}}$,
$\mathbf{K}=\mathbf{P}\mathbf{\Psi}(\mathbf{\Psi}^{T}\mathbf{P}\mathbf{\Psi}+\mathbf{I})^{-1}$
and $\mathbf{\widetilde{\mathbf{K}}}$ is a constant matrix such that
$\mathbf{A}^{T}\left(\mathbf{I}-\widetilde{\mathbf{K}}\mathbf{\Psi}\right)$
is Hurwitz. It follows that there exists 
\begin{align}
 & \mathbf{P}^{\left(2\right)}=\nonumber \\
 & \sum_{i=0}^{\infty}(\mathbf{A}^{T}(\mathbf{I}-\widetilde{\mathbf{K}}\mathbf{\Psi}))^{i}(\mathbf{A}^{T}\widetilde{\mathbf{K}}\mathbf{\widetilde{\mathbf{K}}}^{T}\mathbf{A}+\mathbf{Q})((\mathbf{I}-\widetilde{\mathbf{K}}\mathbf{\Psi})^{T}\mathbf{A})^{i}
\end{align}
 such that $\mathbf{P}^{\left(2\right)}\geq g\left(\mathbf{P}^{\left(2\right)}|\mathbf{H}_{k}=\mathbf{H}\right)>\mathbf{0}$.
Further note that
\begin{align}
\mathbf{\Psi}_{k}\mathbf{\Psi}_{k}^{T} & =\delta_{k}\mathbf{B}\mathbf{H}_{k}\left(\delta_{k}\mathbf{H}_{k}^{T}\mathbf{M}\mathbf{H}_{k}+\mathbf{R}\right)^{-1}\mathbf{H}_{k}^{T}\mathbf{B}^{T}\nonumber \\
 & =\mathbf{B}\left(\mathbf{M}\right)^{-1}\mathbf{B}^{T}-\mathbf{B}\left(\delta_{k}\mathbf{M}\mathbf{H}_{k}\mathbf{R}\mathbf{H}_{k}^{T}\mathbf{M}+\mathbf{M}\right)^{-1}\mathbf{B}^{T}.
\end{align}

As a result, if there exist $\mathbf{H}^{(1)}$ and $\mathbf{H}^{(2)}$
such that $\mathbf{H}^{(1)}\mathbf{R}(\mathbf{H}^{(1)})^{T}>\mathbf{H}^{(2)}\mathbf{R}(\mathbf{H}^{(2)})^{T}$,
we have $\mathbf{\Psi}^{(1)}(\mathbf{\Psi}^{(1)})^{T}>\mathbf{\Psi}^{(2)}(\mathbf{\Psi}^{(2)})^{T}$
and $g\left(\mathbf{P}^{\left(2\right)}|\mathbf{H}_{k}=\mathbf{H}^{(1)}\right)<g\left(\mathbf{P}^{\left(2\right)}|\mathbf{H}_{k}=\mathbf{H}^{(2)}\right)$.
Therefore,
\begin{align}
 & g(\mathbf{P}^{\left(2\right)})=\mathbb{E}[g(\mathbf{P}^{\left(2\right)})|\mathbf{H}_{k}\mathbf{R}\mathbf{H}_{k}^{T}\geq\mathbf{H}\mathbf{R}\mathbf{H}^{T}]+\mathbb{E}[g(\mathbf{P}^{\left(2\right)})|\mathbf{H}_{k}\mathbf{R}\mathbf{H}_{k}^{T}\nonumber \\
 & <\mathbf{H}\mathbf{R}\mathbf{H}^{T}]<\mathbf{P}^{\left(2\right)}\mathrm{Pr}(\mathbf{H}_{k}\mathbf{R}\mathbf{H}_{k}^{T}\geq\mathbf{H}\mathbf{R}\mathbf{H}^{T})+(\mathbf{A}^{T}\mathbf{P}^{\left(2\right)}\mathbf{A}+\mathbf{Q})\nonumber \\
 & \cdot\mathrm{Pr}\left(\mathbf{H}_{k}\mathbf{R}\mathbf{H}_{k}^{T}<\mathbf{H}\mathbf{R}\mathbf{H}^{T}\right).
\end{align}

By letting $\mathbf{H}\rightarrow\mathbf{0}$, we have $\mathrm{Pr}(\mathbf{H}_{k}\mathbf{R}\mathbf{H}_{k}^{T}\geq\mathbf{H}\mathbf{R}\mathbf{H}^{T})\rightarrow1$.
It follows that there exists $\mathbf{P}^{\left(2\right)}\geq g\left(\mathbf{P}^{\left(2\right)}\right)$.
We now construct two matrix sequences: 
\begin{align}
\left\{ \mathbf{P}_{k}^{\left(1\right)}:\mathbf{P}_{k+1}^{\left(1\right)}=g(\mathbf{P}_{k}^{\left(1\right)}),\mathbf{P}_{0}^{\left(1\right)}=\mathbf{P}^{\left(1\right)},k\geq0\right\} ,\label{eq:sequence-1}
\end{align}
\begin{align}
\left\{ \mathbf{P}_{k}^{\left(2\right)}:\mathbf{P}_{k+1}^{\left(2\right)}=g(\mathbf{P}_{k}^{\left(2\right)}),\mathbf{P}_{0}^{\left(2\right)}=\mathbf{P}^{\left(2\right)},k\geq0\right\} .\label{eq:sequence-2}
\end{align}

Due to the monotonicity of $g\left(\cdot\right)$, it follows that
$\mathbf{P}_{k+1}^{\left(1\right)}\geq\mathbf{P}_{k}^{\left(1\right)},\forall k\geq0$,
and $\mathbf{P}_{k+1}^{\left(2\right)}\leq\mathbf{P}_{k}^{\left(2\right)},\forall k\geq0$.
Therefore, we have 
\begin{align}
\mathbf{P}_{k}^{\left(1\right)}\leq\mathbf{P}_{k+1}^{\left(1\right)} & \leq\mathbf{P}_{k+1}^{\left(2\right)}\leq\mathbf{P}_{k}^{\left(2\right)}\leq\mathbf{P}^{\left(2\right)}.
\end{align}

Therefore, the monotonically increasing sequence $\left\{ \mathbf{P}_{k}^{\left(1\right)},k\geq0\right\} $
is bounded from above, i.e., $\mathbf{P}_{k}^{\left(1\right)}\leq\mathbf{P}^{\left(2\right)},\forall k\geq0$,
it follows that the sequence $\left\{ \mathbf{P}_{k}^{\left(1\right)},k\geq0\right\} $
is convergent, i.e., there is a $\left(\mathbf{P}^{\left(1\right)}\right)^{\ast}$
such that
\begin{align}
 & \lim_{k\rightarrow\infty}\mathbf{P}_{k}^{\left(1\right)}=\left(\mathbf{P}^{\left(1\right)}\right)^{\ast}=g\left(\left(\mathbf{P}^{\left(1\right)}\right)^{\ast}\right).
\end{align}
 Therefore, we prove the existence of $\mathbf{P}^{*}$ that satisfies
the NME (\ref{eq: NME}) when one of the three conditions \textsl{(a.1)},
\textsl{(a.2)} and \textsl{(a.3)} in Lemma \ref{lemma-impacts-MIMO Fading}
is satisfied.

We now prove the existence of $\mathbf{P}^{*}$ such that $\mathbf{P}^{*}=g\left(\mathbf{P}^{*}\right)$
under the sufficient condition (\ref{eq:suff-condition-existence}).
Based on Lemma \ref{lemma-decomposed-NME}, we substitute $\boldsymbol{\mathbf{P}}_{k}^{c}$
(\ref{eq:Pk-c}) and $\boldsymbol{\mathbf{P}}_{k}^{uc}$ (\ref{eq:Pk-uc})
into the decomposed NME (\ref{eq:fine-grained-NME}). The $\boldsymbol{\mathbf{P}}_{k}^{uc}$
dependent terms in the decomposed NME (\ref{eq:fine-grained-NME})
can be represented as
\begin{align}
 & \mathbf{A}^{T}\mathbb{E}\left[\boldsymbol{\mathbf{P}}_{k}^{uc}\right]\mathbf{A}=\mathbf{A}^{T}\mathbb{E}[\mathbf{V}_{k}^{T}\left(\mathbf{I}-\boldsymbol{\Pi}_{k}\right)\mathbf{V}_{k}\mathbf{P}\mathbf{V}_{k}^{T}\left(\mathbf{I}-\boldsymbol{\Pi}_{k}\right)\mathbf{V}_{k}]\mathbf{A}\nonumber \\
 & -\mathbf{A}^{T}\mathbb{E}[\mathbf{V}_{k}^{T}\left[\begin{array}{cc}
\mathbf{0}_{\gamma_{k}} & \mathbf{0}_{\gamma_{k}\times\left(S-\gamma_{k}\right)}\\
\mathbf{0}_{\left(S-\gamma_{k}\right)\times\gamma_{k}} & \boldsymbol{\Sigma}_{k}^{T}\left(\mathbf{V}_{k}\mathbf{P}\mathbf{V}_{k}^{H}\right)_{\gamma_{k}}\boldsymbol{\Sigma}_{k}
\end{array}\right]\mathbf{V}_{k}]\mathbf{A}.\label{eq:p-uc-detail}
\end{align}

Moreover, denote $\widetilde{\boldsymbol{\Lambda}_{k}}=\mathrm{diag}\left(\left(\boldsymbol{\Lambda}_{k}\right)_{\gamma_{k}}^{\frac{1}{2}},\mathbf{I}_{\left(S-\gamma_{k}\right)}\right)$,
it follows that 
\begin{align}
 & \delta_{k}\mathbf{B}\mathbf{H}_{k}\left(\delta_{k}\mathbf{H}_{k}^{T}\mathbf{M}\mathbf{H}_{k}+\mathbf{R}\right)^{-1}\mathbf{H}_{k}^{T}\mathbf{B}^{T}=\mathbf{V}_{k}^{T}\widetilde{\boldsymbol{\Lambda}_{k}}\boldsymbol{\Pi}_{k}\widetilde{\boldsymbol{\Lambda}_{k}}\mathbf{V}_{k}.
\end{align}

Denote $\left(\mathbf{V}_{k}\mathbf{P}\mathbf{V}_{k}^{T}\right)_{\gamma_{k}}=\mathbf{P}_{\gamma_{k}}$
and $\left(\boldsymbol{\Lambda}_{k}\right)_{\gamma_{k}}^{\frac{1}{2}}\left(\mathbf{V}_{k}\mathbf{P}\mathbf{V}_{k}^{T}\right)_{\gamma_{k}}\left(\boldsymbol{\Lambda}_{k}\right)_{\gamma_{k}}^{\frac{1}{2}}=\widetilde{\mathbf{P}}_{\gamma_{k}}.$
The $\boldsymbol{\mathbf{P}}_{k}^{c}$ dependent terms in the decomposed
NME (\ref{eq:fine-grained-NME}) thus can be represented as
\begin{align}
 & \mathbf{A}^{T}\mathbb{E}[\boldsymbol{\mathbf{P}}_{k}^{c}\mathbf{A}-\delta_{k}\mathbf{A}^{T}\boldsymbol{\mathbf{P}}_{k}^{c}\mathbf{B}\mathbf{H}_{k}(\delta_{k}\mathbf{H}_{k}^{T}\mathbf{B}^{T}\boldsymbol{\mathbf{P}}_{k}^{c}\mathbf{B}\mathbf{H}_{k}\nonumber \\
 & +\delta_{k}\mathbf{H}_{k}^{T}\mathbf{M}\mathbf{H}_{k}+\mathbf{R})^{-1}\mathbf{H}_{k}^{T}\mathbf{B}^{T}\boldsymbol{\mathbf{P}}_{k}^{c}]\mathbf{A}\nonumber \\
 & =\mathbf{A}^{T}\mathbb{E}[(\boldsymbol{\mathbf{P}}_{k}^{c}\mathbf{\Psi}_{k}\mathbf{\Psi}_{k}^{T}+\mathbf{I})^{-1}\boldsymbol{\mathbf{P}}_{k}^{c}]\mathbf{A}\nonumber \\
 & =\mathbf{A}^{T}\mathbb{E}[\mathbf{V}_{k}^{T}(\mathbf{V}_{k}\boldsymbol{\mathbf{P}}_{k}^{c}\mathbf{V}_{k}^{T}\boldsymbol{\Lambda}_{k}+\mathbf{I})^{-1}\mathbf{V}_{k}\boldsymbol{\mathbf{P}}_{k}^{c}\mathbf{V}_{k}^{T}\mathbf{V}_{k}]\mathbf{A}\nonumber \\
 & =\mathbf{A}^{T}\mathbb{E}[\mathbf{V}_{k}^{T}\widetilde{\boldsymbol{\Lambda}_{k}}^{-1}(\widetilde{\boldsymbol{\Lambda}_{k}}\mathbf{V}_{k}\boldsymbol{\mathbf{P}}_{k}^{c}\mathbf{V}_{k}^{T}\widetilde{\boldsymbol{\Lambda}_{k}}\boldsymbol{\Pi}_{k}+\mathbf{I})^{-1}\nonumber \\
 & \cdot\widetilde{\boldsymbol{\Lambda}_{k}}\mathbf{V}_{k}\boldsymbol{\mathbf{P}}_{k}^{c}\mathbf{V}_{k}^{T}\widetilde{\boldsymbol{\Lambda}_{k}}\widetilde{\boldsymbol{\Lambda}_{k}}^{-1}\mathbf{V}_{k}]\mathbf{A}=\mathbf{A}^{T}\mathbb{E}[\mathbf{V}_{k}^{T}\widetilde{\boldsymbol{\Lambda}_{k}}^{-1}\nonumber \\
 & \left[\begin{array}{cc}
\left(\mathbf{I}+\widetilde{\mathbf{P}}_{\gamma_{k}}^{-1}\right)^{-1} & \left(\mathbf{I}+\widetilde{\mathbf{P}}_{\gamma_{k}}^{-1}\right)^{-1}\mathbf{P}_{\gamma_{k}}\boldsymbol{\Sigma}_{k}\\
\boldsymbol{\Sigma}_{k}^{T}\mathbf{P}_{\gamma_{k}}\left(\mathbf{I}+\widetilde{\mathbf{P}}_{\gamma_{k}}^{-1}\right)^{-1} & \mathbf{\Theta}_{k}
\end{array}\right]\nonumber \\
 & \widetilde{\boldsymbol{\Lambda}_{k}}^{-1}\mathbf{V}_{k}]\mathbf{A},\label{eq:p-c-detail}
\end{align}
where 
\begin{align}
 & \mathbf{\Theta}_{k}=\boldsymbol{\Sigma}_{k}^{T}\mathbf{P}_{\gamma_{k}}\boldsymbol{\Sigma}_{k}-\boldsymbol{\Sigma}_{k}^{T}\mathbf{P}_{\gamma_{k}}\left(\mathbf{I}+\widetilde{\mathbf{P}}_{\gamma_{k}}^{-1}\right)^{-1}\mathbf{P}_{\gamma_{k}}\boldsymbol{\Sigma}_{k}.
\end{align}

Substituting (\ref{eq:p-uc-detail}) and (\ref{eq:p-c-detail}) into
(\ref{eq:fine-grained-NME}), it follows that
\begin{align}
 & g(\mathbf{P})=\mathbf{Q}+\nonumber \\
 & \mathbf{A}^{T}\mathbb{E}[\mathbf{V}_{k}^{T}\left(\mathbf{I}-\boldsymbol{\Pi}_{k}\right)\mathbf{V}_{k}\mathbf{P}\mathbf{V}_{k}^{T}\left(\mathbf{I}-\boldsymbol{\Pi}_{k}\right)\mathbf{V}_{k}]\mathbf{A}+\mathbf{A}^{T}\mathbb{E}[\mathbf{V}_{k}^{T}\nonumber \\
 & \cdot\left[\begin{array}{cc}
\left(\boldsymbol{\Lambda}_{k}\right)_{\gamma_{k}}^{-\frac{1}{2}}\left(\mathbf{I}+\mathbf{P}_{\gamma_{k}}^{-1}\right)^{-1}\left(\boldsymbol{\Lambda}_{k}\right)_{\gamma_{k}}^{-\frac{1}{2}} & \left(\mathbf{I}+\widetilde{\mathbf{P}}_{\gamma_{k}}^{-1}\right)^{-1}\mathbf{P}_{\gamma_{k}}\boldsymbol{\Sigma}_{k}\\
\boldsymbol{\Sigma}_{k}^{T}\mathbf{P}_{\gamma_{k}}\left(\mathbf{I}+\widetilde{\mathbf{P}}_{\gamma_{k}}^{-1}\right)^{-1} & -\mathbf{\widetilde{\Theta}}_{k}
\end{array}\right]\nonumber \\
 & \cdot\mathbf{V}_{k}]\mathbf{A},
\end{align}
where 
\begin{align}
 & \mathbf{\widetilde{\Theta}}_{k}=\boldsymbol{\Sigma}_{k}^{T}\mathbf{P}_{\gamma_{k}}\left(\mathbf{I}+\widetilde{\mathbf{P}}_{\gamma_{k}}^{-1}\right)^{-1}\mathbf{P}_{\gamma_{k}}\boldsymbol{\Sigma}_{k}\geq\mathbf{0}.
\end{align}

Note that $\left(\mathbf{I}+\mathbf{P}_{\gamma_{k}}^{-1}\right)^{-1}\leq\mathbf{I}$,
it follows that
\begin{align}
g(\mathbf{P})< & \mathbf{Q}+\mathbf{A}^{T}\mathbb{E}[\mathbf{V}_{k}^{T}\left(\mathbf{I}-\boldsymbol{\Pi}_{k}\right)\mathbf{V}_{k}\mathbf{P}\nonumber \\
 & \cdot\mathbf{V}_{k}^{T}\left(\mathbf{I}-\boldsymbol{\Pi}_{k}\right)\mathbf{V}_{k}]\mathbf{A}+\left\Vert \mathbf{A}\right\Vert ^{2}\mathbb{E}[\mathrm{Tr}(\left(\boldsymbol{\Lambda}_{k}\right)_{\gamma_{k}}^{-1})]\mathbf{I}
\end{align}

Therefore, under condition (\ref{eq:suff-condition-existence}), there
is a $\mathbf{P}^{\left(2\right)}=\vartheta\mathbf{I}$ with 
\begin{align}
 & \vartheta\geq\frac{\left\Vert \mathbf{Q}\right\Vert +\left\Vert \mathbf{A}\right\Vert ^{2}\mathbb{E}[\mathrm{Tr}(\left(\boldsymbol{\Lambda}_{k}\right)_{\gamma_{k}}^{-1})]}{1-\left\Vert \mathbb{E}\left[\mathbf{A}^{T}\mathbf{V}_{k}^{T}\left(\mathbf{I}-\boldsymbol{\Pi}_{k}\right)\mathbf{V}_{k}\mathbf{A}\right]\right\Vert },
\end{align}
such that $g\left(\mathbf{P}^{\left(2\right)}\right)\leq\mathbf{P}^{\left(2\right)}$.
Using similar techniques and constructing the two matrix sequences
the same way as in (\ref{eq:sequence-1}) and (\ref{eq:sequence-2}),
it follows that 
\begin{align}
 & \lim_{k\rightarrow\infty}\mathbf{P}_{k}^{\left(1\right)}=\left(\mathbf{P}^{\left(1\right)}\right)^{\ast}=g\left(\left(\mathbf{P}^{\left(1\right)}\right)^{\ast}\right).
\end{align}
 Therefore, the existence of $\mathbf{P}^{*}$ that satisfies the
NME (\ref{eq: NME}) under the sufficient condition (\ref{eq:suff-condition-existence})
in Lemma \ref{lemma-impacts-MIMO Fading} is proved.

In the following, we shall prove the uniqueness of $\mathbf{P}^{*}$.
Suppose there are $\left(\mathbf{P}^{\left(1\right)}\right)^{\ast}$
and $\left(\mathbf{P}^{\left(2\right)}\right)^{\ast}$ such that $\left(\mathbf{P}^{\left(1\right)}\right)^{\ast}=g\left(\left(\mathbf{P}^{\left(1\right)}\right)^{\ast}\right)$
and $\left(\mathbf{P}^{\left(2\right)}\right)^{\ast}=g\left(\left(\mathbf{P}^{\left(2\right)}\right)^{\ast}\right)$.
Then, there is a positive constant $\phi^{*}\in(0,1)$ such that $\left(\mathbf{P}^{\left(1\right)}\right)^{\ast}\geq\phi^{*}\left(\mathbf{P}^{\left(2\right)}\right)^{\ast}$
and $\left(\mathbf{P}^{\left(1\right)}\right)^{\ast}\ngeq\phi\left(\mathbf{P}^{\left(2\right)}\right)^{\ast}$
for all $\phi>\phi^{*}$. Note that 
\begin{align}
 & g(\phi^{*}(\mathbf{P}^{\left(2\right)})^{\ast})=\mathbb{E}[\mathbf{A}(\mathbf{\Psi}_{k}\mathbf{\Psi}_{k}^{T}+(\phi^{*}(\mathbf{P}^{\left(2\right)})^{\ast})^{-1})^{-1}\mathbf{A}^{T}]+\mathbf{Q}\nonumber \\
 & \geq\mathbb{E}[(\phi^{*})^{-1}\mathbf{A}(\mathbf{\Psi}_{k}\mathbf{\Psi}_{k}^{T}+(\phi^{*}(\mathbf{P}^{\left(2\right)})^{\ast})^{-1})^{-1}\mathbf{A}^{T}]+\mathbf{Q}\nonumber \\
 & \geq\left(1+\varphi\right)\phi^{*}g((\mathbf{P}^{\left(2\right)})^{\ast}),
\end{align}
where $\varphi=\frac{(1-\phi^{*})\sigma_{\mathbf{Q}}}{|||g((\mathbf{P}^{\left(2\right)})^{\ast})||}$
is a positive constant and $\sigma_{\mathbf{Q}}$ is the minimum singular
value of $\mathbf{Q}$. Further note that 
\begin{align}
 & (\mathbf{P}^{\left(1\right)})^{\ast}=g((\mathbf{P}^{\left(1\right)})^{\ast})\geq g(\phi^{*}(\mathbf{P}^{\left(2\right)})^{\ast})\nonumber \\
 & \geq\left(1+\varphi\right)\phi^{*}g((\mathbf{P}^{\left(2\right)})^{\ast})=\left(1+\varphi\right)\phi^{*}\left(\mathbf{P}^{\left(2\right)}\right)^{\ast}.
\end{align}
This means that there is a $\phi=\left(1+\varphi\right)\phi^{*}>\phi^{*}$
such that $\left(\mathbf{P}^{\left(1\right)}\right)^{\ast}\geq\phi\left(\mathbf{P}^{\left(2\right)}\right)^{\ast}$,
which contradicts the fact that $\left(\mathbf{P}^{\left(1\right)}\right)^{\ast}\ngeq\phi\left(\mathbf{P}^{\left(2\right)}\right)^{\ast}$
for all $\phi>\phi^{*}$. As a result, the uniqueness of $\mathbf{P}^{*}$
that satisfies $\left(\mathbf{P}\right)^{\ast}=g\left(\left(\mathbf{P}\right)^{\ast}\right)$
is proved. Therefore, Theorem \ref{Thm: suff-condition-existence-uniqueness}
is proved.

\subsection{\label{subsec:Proof-of-Lemma-4-Lemma-5}Proof of Lemma \ref{Lemma: properties-iter}
and Lemma \ref{Lemma: ODE Asympt-Stability}}

Note that
\begin{align}
 & f(\mathbf{P}^{(1)})-f(\mathbf{P}^{(2)})=\mathbf{P}^{(2)}-\mathbf{P}^{(1)}\nonumber \\
 & \mathbb{E}[\mathbf{A}^{T}(\mathbf{I}-\mathbf{K}_{k}^{(1)}\mathbf{\Psi}_{k})\mathbf{P}^{(1)}(\mathbf{I}-\mathbf{K}_{k}^{(1)}\mathbf{\Psi}_{k})\mathbf{A}+\mathbf{K}_{k}^{(1)}(\mathbf{K}_{k}^{(1)})^{T}]-\nonumber \\
 & \mathbb{E}[\mathbf{A}^{T}(\mathbf{I}-\mathbf{K}_{k}^{(2)}\mathbf{\Psi}_{k})\mathbf{P}^{(2)}(\mathbf{I}-\mathbf{K}_{k}^{(2)}\mathbf{\Psi}_{k})^{T}\mathbf{A}+\mathbf{K}_{k}^{(2)}(\mathbf{K}_{k}^{(2)})^{T}]\nonumber \\
 & \leq\mathbb{E}[\mathbf{A}^{T}(\mathbf{I}-\mathbf{K}_{k}^{(2)}\mathbf{\Psi}_{k})(\mathbf{P}^{(1)}-\mathbf{P}^{(2)})(\mathbf{I}-\mathbf{K}_{k}^{(2)}\mathbf{\Psi}_{k})^{T}\mathbf{A}].
\end{align}
Note that $\left\Vert \mathbf{A}^{T}\left(\mathbf{I}-\mathbf{K}_{k}^{(2)}\mathbf{\Psi}_{k}\right)\right\Vert \leq\left\Vert \mathbf{A}\right\Vert $,
it follows that 
\begin{align}
 & \left\Vert f(\mathbf{P}^{(1)})-f(\mathbf{P}^{(2)})\right\Vert \nonumber \\
 & \leq\left\Vert \mathbf{P}^{(1)}-\mathbf{P}^{(2)}\right\Vert +\mathbb{E}\left[\left\Vert \mathbf{A}^{T}\left(\mathbf{I}-\mathbf{K}_{k}^{(2)}\mathbf{\Psi}_{k}\right)\right\Vert ^{2}\left\Vert \mathbf{P}^{(1)}-\mathbf{P}^{(2)}\right\Vert \right]\nonumber \\
 & \leq\left(1+\left\Vert \mathbf{A}\right\Vert ^{2}\right)\left\Vert \mathbf{P}^{(1)}-\mathbf{P}^{(2)}\right\Vert .
\end{align}

Note that for any given realization of $\mathbf{P}_{k}$, $\mathbf{N}_{k+1}=\left(\widehat{f}\left(\mathbf{P}_{k+1}\right)-f\left(\mathbf{P}_{k+1}\right)\right)$
is a function of $\left\{ \delta_{k+1},\mathbf{H}_{k+1}\right\} $.
Moreover, due to the i.i.d. property of $\left\{ \delta_{k+1},\mathbf{H}_{k+1}\right\} $,
it follows that $\mathbb{E}\left[\left.\widehat{f}\left(\mathbf{P}_{k+1}\right)-f\left(\mathbf{P}_{k+1}\right)\right|\mathbf{P}_{k}\right]=\mathbf{0}$.
Therefore, $\mathbb{E}\left[\left.\mathbf{N}_{k+1}\right|\mathcal{F}_{k}\right]=\mathbf{0},\forall k>0.$

Note that 
\begin{align}
 & \mathbb{E}[||\mathbf{N}_{k+1}||^{2}|\mathbf{P}_{k}]\leq\mathbb{E}[||\widehat{f}\left(\mathbf{P}_{k+1}\right)-\mathbf{Q}||^{2}+||f\left(\mathbf{P}_{k+1}\right)\nonumber \\
 & -\mathbf{Q}||^{2}|\mathbf{P}_{k}]\leq\mathbb{E}[||\mathbf{A}(\mathbf{\Psi}_{k}\mathbf{\Psi}_{k}^{T}+(\mathbf{P}_{k})^{-1})^{-1}\mathbf{A}^{T}||^{2}|\mathbf{P}_{k}]\nonumber \\
 & +||\mathbb{E}[\mathbf{A}(\mathbf{\Psi}_{k}\mathbf{\Psi}_{k}^{T}+(\mathbf{P}_{k})^{-1})^{-1}\mathbf{A}^{T}|\mathbf{P}_{k}]||^{2}\nonumber \\
 & \leq2||\mathbf{A}\mathbf{P}_{k}\mathbf{A}^{T}||^{2}\leq2\left\Vert \mathbf{A}\right\Vert ^{2}\left\Vert \mathbf{P}_{k}\right\Vert ^{2}.
\end{align}
Therefore, Lemma \ref{Lemma: properties-iter} is proved.

According to Theorem 2.1 in Chapter 5 of {[}38{]}, when all the conditions
in Lemma \ref{Lemma: properties-iter} are satisfied, if the limiting
ODE (\ref{eq: SA-ODE}) has a unique equilibrium point $\mathbf{P}^{*}$
that is globally asymptotically stable, then $\mathbf{P}_{k}$ converges
to $\mathbf{P}^{*}$ with probability 1. Therefore, Lemma \ref{Lemma: ODE Asympt-Stability}
is proved.

\subsection{\label{subsec:Proof-of-Lemma-6}Proof of Lemma \ref{Lemma: Difference State trajectory}}

Let $L=N\xi$ for some $N>0$. For $t>0$, denote $\left[t\right]=\max\left\{ k\xi:n>0,k\xi<t\right\} $.
For $n\geq0$ and $1\leq l\leq L$, we have 
\begin{align}
 & \overline{\mathbf{P}}\left(t_{k+l}\right)=\overline{\mathbf{P}}\left(t_{k}\right)+\intop_{t_{k}}^{t_{k+l}}f\left(\overline{\mathbf{P}}\left(\left[t\right]\right)\right)\mathrm{d}t,\label{eq:lemma-6-1}\\
 & \mathbf{P}^{t_{k}}\left(t_{k+l}\right)=\overline{\mathbf{P}}\left(t_{k}\right)+\intop_{t_{k}}^{t_{k+l}}f\left(\mathbf{P}^{t_{k}}\left(\left[t\right]\right)\right)\mathrm{d}t\nonumber \\
 & +\intop_{t_{k}}^{t_{k+l}}\left(f\left(\mathbf{P}^{t_{k}}\left(t\right)\right)-f\left(\mathbf{P}^{t_{k}}\left(\left[t\right]\right)\right)\right)\mathrm{d}t.\label{eq:lemma-6-2}
\end{align}
Subtracting (\ref{eq:lemma-6-2}) from (\ref{eq:lemma-6-1}) and noting
that 
\begin{align}
 & ||\intop_{t_{k}}^{t_{k+l}}\left(f\left(\mathbf{P}^{t_{k}}\left(t\right)\right)-f\left(\mathbf{P}^{t_{k}}\left(\left[t\right]\right)\right)\right)\mathrm{d}t||\nonumber \\
 & \leq\left\Vert \mathbf{A}\right\Vert ^{2}\xi\sum_{m=0}^{l-1}\sup_{j\leq m}||\overline{\mathbf{P}}\left(t_{k+j}\right)-\mathbf{P}^{t_{k}}\left(t_{k+j}\right)||,\\
 & ||\intop_{t_{k+l}}^{t_{k+l+1}}\left(f\left(\mathbf{P}\left(t\right)\right)-f\left(\mathbf{P}\left(\left[t\right]\right)\right)\right)\mathrm{d}t||\leq c_{1}\xi\left(1+\overline{\mathbf{P}}\left(t_{k}\right)\right),
\end{align}
it follows that 
\begin{align}
 & \sup_{0\leq j\leq l}||\overline{\mathbf{P}}\left(t_{k+j}\right)-\mathbf{P}^{t_{k}}\left(t_{k+j}\right)||\leq c_{1}\xi\left(1+\overline{\mathbf{P}}\left(t_{k}\right)\right)\nonumber \\
 & +\xi L\left\Vert \mathbf{A}\right\Vert ^{2}\sum_{m=0}^{l-1}\sup_{j\leq m}||\overline{\mathbf{P}}\left(t_{k+j}\right)-\mathbf{P}^{t_{k}}\left(t_{k+j}\right)||.
\end{align}
By the discrete Gronwall inequality, it follows that 
\begin{align}
 & \sup_{k\leq j\leq k+N}||\overline{\mathbf{P}}\left(t_{j}\right)-\mathbf{P}^{t_{k}}\left(t_{j}\right)||^{2}\leq c_{2}\xi,
\end{align}
where $c_{2}$ is a constant. Since both $\sup_{t_{j}\leq t\leq t_{j+1}}||\overline{\mathbf{P}}\left(t\right)-\overline{\mathbf{P}}\left(t_{j}\right)||^{2}$
and $\sup_{t_{j}\leq t\leq t_{j+1}}||\mathbf{P}^{t_{k}}\left(t\right)-\mathbf{P}^{t_{k}}\left(t_{j}\right)||^{2}$
are $\mathcal{O}\left(\xi\right)$, it follows that
\begin{align}
 & \sup_{t\in\left[0,L\right]}\left\Vert \overline{\mathbf{P}}\left(l+t\right)-\mathbf{P}^{l}\left(l+t\right)\right\Vert \leq c_{3}\xi,
\end{align}
where $c_{3}$ is a constant. Therefore, Lemma \ref{Lemma: Difference State trajectory}
is proved.

\subsection{\label{subsec:Proof-of-Thm-6}Proof of Theorem \ref{Thm: Convergence results}}

Note that if the virtual fixed-point process $\left\{ \widetilde{\mathbf{P}}_{k},k\geq0\right\} $
in (\ref{eq: virtual-fixed-point-process}) corresponds to the fixed-point
equation $\mathbf{P}=\widetilde{g}(\mathbf{P})$ with $\widetilde{g}(\mathbf{P})=\mathbf{P}+\xi f\left(\mathbf{P}\right)$,
we know that if one of the three conditions \textsl{(a.1)}, \textsl{(a.2)}
and \textsl{(a.3)} in Lemma \ref{lemma-impacts-MIMO Fading} is satisfied,
or the condition (\ref{eq:suff-condition-existence}) in Theorem \ref{Thm: suff-condition-existence-uniqueness}
is satisfied, the solution $\mathbf{P}^{*}$ to the fixed-point equation
$\mathbf{P}^{*}=\widetilde{g}(\mathbf{P}^{*})$ exists and is unique.
Using similar techniques as in Appendix \ref{subsec:Proof-of-Thm-5},
there is a $\widetilde{\mathbf{P}}^{\left(1\right)}=\mathbf{0}$ such
that $\widetilde{\mathbf{P}}^{\left(1\right)}<\widetilde{g}\left(\widetilde{\mathbf{P}}^{\left(1\right)}\right)$,
and a sufficiently large $\widetilde{\mathbf{P}}^{\left(2\right)}$
such that $\widetilde{\mathbf{P}}^{\left(2\right)}>\widetilde{g}\left(\widetilde{\mathbf{P}}^{\left(2\right)}\right)$.
We now construct the following two matrix sequences: 
\begin{align}
 & \left\{ \widetilde{\mathbf{P}}_{k}^{\left(1\right)}:\widetilde{\mathbf{P}}_{k+1}^{\left(1\right)}=\widetilde{g}(\widetilde{\mathbf{P}}_{k}^{\left(1\right)}),\mathbf{P}_{0}^{\left(1\right)}=\mathbf{0},k\geq0\right\} ,\label{eq:sequence-1-1}\\
 & \left\{ \widetilde{\mathbf{P}}_{k}^{\left(2\right)}:\widetilde{\mathbf{P}}_{k+1}^{\left(2\right)}=\widetilde{g}(\widetilde{\mathbf{P}}_{k}^{\left(2\right)}),\mathbf{P}_{0}^{\left(2\right)}=\widetilde{\mathbf{P}}^{\left(2\right)},k\geq0\right\} .\label{eq:sequence-2-1}
\end{align}

Let the initial condition of the fixed-point process be $\mathbf{0}\leq\widetilde{\mathbf{P}}_{0}\leq\widetilde{\mathbf{P}}^{\left(2\right)}$,
it follows that $\widetilde{\mathbf{P}}_{k}^{\left(1\right)}\leq\widetilde{\mathbf{P}}_{k}\leq\widetilde{\mathbf{P}}_{k}^{\left(2\right)}$.
Let $k\rightarrow\infty$ and note that $\mathbf{P}^{*}$ exists and
is unique, it follows that
\begin{align}
 & \mathbf{P}^{*}=\lim_{k\rightarrow\infty}\widetilde{\mathbf{P}}_{k}^{\left(1\right)}\leq\lim_{k\rightarrow\infty}\widetilde{\mathbf{P}}_{k}\leq\lim_{k\rightarrow\infty}\widetilde{\mathbf{P}}_{k}^{\left(2\right)}=\mathbf{P}^{*}.
\end{align}

Since $\widetilde{\mathbf{P}}^{\left(2\right)}$ can be arbitrarily
large, it follows that for any bounded initial value $\widetilde{\mathbf{P}}_{0}$,
the virtual fixed-point process $\left\{ \widetilde{\mathbf{P}}_{k},k\geq0\right\} $
in (\ref{eq: virtual-fixed-point-process}) converges to $\mathbf{P}^{*}$.
Based on Lemma \ref{Lemma: Difference State trajectory}, the limiting
ODE (\ref{eq: SA-ODE}) thus has a unique equilibrium point $\mathbf{P}^{*}$
that is globally asymptotically stable. Moreover, based on Lemma \ref{Lemma: ODE Asympt-Stability},
it follows that the $\mathbf{P}_{k}$ obtained by stochastic approximation
iteration (\ref{eq: prps-SA}) converges to $\mathbf{P}^{*}$ almost
surely. Based on the structural properties (\ref{eq:opt value function})
and (\ref{eq:optimal control action}) in Theorem \ref{Thm: structure-property-reduced-state-Bellman},
it follows that $\widetilde{V}_{k}\left(\mathbf{x}_{k}\right)$ and
$\mathbf{u}_{k}^{*}\left(\mathbf{\mathbf{S}}_{k}\right)$ converges
to the optimal value function $\widetilde{V}\left(\mathbf{x}_{k}\right)$
and optimal control action $\mathbf{u}^{*}\left(\mathbf{\mathbf{S}}_{k}\right)$
w.p.1., respectively. Therefore, Theorem \ref{Thm: Convergence results}
is proved.

\end{document}